\documentclass[final]{siamltex}
\usepackage{amsmath, amssymb, amsfonts}
\usepackage{graphicx,color,caption,subcaption}
\usepackage{diagbox}
\usepackage[ruled]{algorithm2e}
\usepackage{algorithmic}

\usepackage{enumitem}
\usepackage{multirow}

\newtheorem{remark}{Remark}

\newcommand{\bvec}[1]{\mathbf{#1}}

\newcommand{\mc}[1]{\mathcal{#1}}
\newcommand{\mcV}{\mathcal{V}}

\newcommand{\xc}{\mathrm{xc}}
\newcommand{\vr}{\bvec{r}}
\newcommand{\vF}{\bvec{F}}

\newcommand{\vR}{\bvec{R}}

\newcommand{\ud}{\,\mathrm{d}}

\newcommand{\KS}{\mathrm{KS}}

\newcommand{\abs}[1]{\lvert#1\rvert}
\newcommand{\norm}[1]{\lVert#1\rVert}

\newcommand{\wt}[1]{\widetilde{#1}}

\newcommand{\ie}{\textit{i.e.}}

\newcommand{\Or}{\mathcal{O}}

\newcommand{\I}{\imath}

\title{Adaptively compressed polarizability operator \\for accelerating
large scale \\ 
\textit{ab initio} phonon calculations}

\author{
	Lin Lin\thanks{Department of Mathematics, University of California, Berkeley, Berkeley, CA 94720 and Computational Research Division, Lawrence Berkeley National Laboratory, Berkeley, CA 94720. Email: \texttt{linlin@math.berkeley.edu}}
\and Ze Xu\thanks{Department of Mathematics, University of California, Berkeley, Berkeley, CA 94720. Email: \texttt{zexu@math.berkeley.edu}}\and
	Lexing Ying\thanks{Department of Mathematics and Institute for Computational and Mathematical Engineering, Stanford University, Stanford, CA 94305. Email: \texttt{lexing@math.stanford.edu}}
}

\begin{document}

\maketitle

\begin{abstract}
Phonon calculations based on first principle electronic structure theory, such as the Kohn-Sham density functional theory, have wide applications in physics, chemistry and material science. The computational cost of first principle phonon calculations  typically scales steeply as $\Or(N_e^4)$, where $N_e$ is the number of electrons in the system. In this work, we develop a new method to reduce the computational complexity of computing the full dynamical matrix, and hence the phonon spectrum, to $\Or(N_e^3)$. The key concept for achieving this is to compress the polarizability operator adaptively with respect to the perturbation of the potential due to the change of the atomic configuration. Such adaptively compressed polarizability operator (ACP) allows accurate computation of the phonon spectrum. The reduction of complexity only weakly depends on the size of the band gap, and our method is applicable to insulators as well as semiconductors with small band gaps. We demonstrate the effectiveness of our method using one-dimensional and two-dimensional model problems.


\end{abstract}

\begin{keywords} 
Density functional perturbation theory, phonon calculations, adaptive compression, polarizability operator, Steinheimer equation, Dyson equation.
\end{keywords}

\begin{AMS}
65F10,65F30,65Z05
\end{AMS} 

\pagestyle{myheadings}
\thispagestyle{plain}
\markboth{L. Lin, Z. Xu and L. Ying}{Adaptively Compressed Polarizability Operator for Large Scale Phonon Calculations}

\section{Introduction}\label{sec:intro}

Kohn-Sham density functional theory
(KSDFT)~\cite{HohenbergKohn1964,KohnSham1965} is the most widely used
electronic structure theory for molecules and systems in condensed
phase. In principle, KSDFT provides an exact description of ground state properties of a many body quantum system, such as electron density, energy, and atomic forces. 
Once the electronic ground state is obtained, many physical and chemical properties of the system can be described by studying the \textit{response} of the quantum system under small perturbation. The theory for describing such response behavior is called the density functional perturbation theory (DFPT)~\cite{BaroniGiannozziTesta1987,GonzeLee1997,BaroniGironcoliDalEtAl2001}.

One important application of DFPT is the description of lattice vibrations.  In the Born-Oppenheimer approximation, lattice vibrations can be described by the dynamical matrix, which is related to the Hessian matrix of the ground state energy with respect to the atomic positions. The eigenfunctions of the dynamical matrix give the phonon modes, and the eigenvalues give the phonon frequencies. 
A large variety of physical properties of solids depend on such phonon calculations. A few examples include infrared spectroscopy, elastic neutron scattering, specific heat, heat conduction, and electron-phonon interaction related behaviors such as superconductivity~\cite{GonzeLee1997,BaroniGironcoliDalEtAl2001}.
Furthermore, the computational procedure of phonon calculations are largely transferable to the calculation of other types of response behavior, such as response to homogeneous electric fields, piezoelectric properties, magnons, and many body perturbation theory for the description of electrons at excited states such as the GW theory~\cite{Hedin1965,OnidaReiningRubio2002}.

Mathematically, the procedure for phonon calculations can be straightforward. When atoms are at their equilibrium positions, the atomic forces (i.e. first order derivatives of the energy with respect to the atomic position) are zero for all atoms. 
To compute the Hessian matrix, one can move one atom at a time slightly
away from its equilibrium position, and compute the corresponding atomic
forces.  This amounts to the finite difference (FD) approximation of the
Hessian matrix, and is referred to as the ``frozen
phonon''~\cite{YinCohen1982,BaroniResta1986} approach in physics. The
FD approach is simple to implement and can be used to
obtain phonon spectrum quickly for systems of small sizes. However,
FD requires in total $d\times N_A\sim \Or(N_e)$ KSDFT calculations, where
$d$ is the spatial dimension (usually $d=3$), and $N_A$ is the number of
atoms. The computational complexity of a single KSDFT calculation
typically scales as $\Or(N_e^3)$, where $N_e$ is the number of electrons
in the system. Since $N_A\sim \Or(N_e)$,  the total cost of the FD
approximation is $\Or(N_e^4)$. This is prohibitively expensive for
systems of large sizes.  Furthermore, the accuracy of the FD approximation is limited by the size of the perturbation,
which cannot be too small due to the numerical noise in the evaluation
of the atomic forces in KSDFT calculations (usually the accuracy of
forces is set to be $10^{-4}\sim 10^{-3}$ Hartree/Bohr). Such
numerical noise also makes it difficult to compute non-linear response
properties, which can require even higher order derivatives of the
energy.

DFPT, on the other hand, can be viewed as the ``proper'' way for
computing derivative quantities in the context of KSDFT. The central
quantity in DFPT is the polarizability operator, which characterizes the
linear response of the electron density with respect to the perturbation
of the external potential. More specifically, phonon calculations
requires applying the polarizability operator to $d\times N_A$ perturbation vectors induced by the change of the atomic configuration.
The polarizability operator can be obtained by solving a Dyson equation iteratively~\cite{BaroniGironcoliDalEtAl2001}, and each iteration step requires the solutions to
$\Or(N_e^2)$ Steinheimer equations.  In general the complexity of DFPT is still $\Or(N_e^4)$. So the main advantage of DFPT is that it gives accurate linear response properties. Furthermore, the same framework can be used to compute non-linear response properties~\cite{GonzeVigneron1989,Gonze1995,BaroniGironcoliDalEtAl2001}. The mathematical aspect of DFPT for reduced Hartree-Fock model systems have recently been analyzed~\cite{CancesMourad2014}. It is also possible to reduce the computational complexity of phonon calculations by ``linear scaling methods''~\cite{Goedecker1999,BowlerMiyazaki2012}. Such methods can be successful in reducing the computational cost for systems of large sizes with substantial band gaps, but this can be challenging for medium-sized systems with relatively small band gaps.

The main computational bottleneck of DFPT is the solution of the $\Or(N_e^2)$ Steinheimer equations~\cite{BaroniGironcoliDalEtAl2001}. KSDFT can be defined as a nonlinear eigenvalue problem with $\Or(N_e)$ eigenfunctions. Each of the $\Or(N_e^2)$ equations in DFPT represents the response of an eigenfunction to a different external perturbation. Hence at first sight it is not possible to reduce the number of equations. However, as $N_e$ becomes large, there will be asymptotically more equations to solve than the size of the matrix.
Hence there is potential room to obtain a set of ``compressed perturbations'', which leads to methods for solving DFPT with lower complexity.

At this point it might be enticing to compress the $\Or(N_e^2)$ equations using standard compression schemes such as singular value decomposition (SVD). However, there is some immediate difficulty associated with SVD type of compression scheme:  The matrix to be compressed is of size $\Or(N_e^2)\times \Or(N_e)$ and of approximate rank $\Or(N_e)$. The associated cost of the SVD type of compression is $\Or(N_e^4)$, and hence there is no saving in asymptotic complexity. 
Furthermore, the $\Or(N_e^2)$ equations need to be solved self-consistently according to the Dyson equation. Hence the initially compressed vectors might not be applicable anymore as the iteration proceeds towards the converged solution. This leads to inaccurate phonon calculations.

In this paper, we develop a new method called the adaptively compressed polarizability operator (ACP) to overcome the above difficulties.  
ACP reduces the complexity for applying the polarizability operator to $\Or(N_e)$ vectors as follows. 1) ACP compresses the $\Or(N_e^2)$ right hand side vectors of the Steinheimer equations into $\Or(N_e)$ vectors, using a recently developed 
randomized density fitting method~\cite{LuYing2015}. Together with a
Chebyshev interpolation procedure to disentangle the energy dependence
from the right hand side vectors, ACP reduces the number of equations
from $\Or(N_e^2)$ to only $\Or(N_e)$. 2) ACP reformulates the Dyson
equation into an equivalent fixed point problem, where the compression of the
Steinheimer equations depends adaptively on the unknown solutions. Using
such adaptive compression procedure, we demonstrate that the
self-consistent solution to the Dyson equation no longer hinders the accuracy of the compressed polarizability operator. Such adaptive compression
strategy shares similar spirit to the recently developed adaptively
compressed exchange operator (ACE) for accelerating KSDFT calculations
with hybrid exchange-correlation functionals~\cite{Lin2016ACE}. We
demonstrate that the overall computational complexity for phonon
calculations can be reduced to $\Or(N_e^3)$, and the cost depends only
weakly on the band gap of the system. Hence the method can be applied to
both insulators and semiconductors with small gaps.  To the extent of
our knowledge, this is the first result of this type in literature.
We demonstrate the numerical performance of the ACP formulation for
accelerating phonon calculations using model systems for one-dimensional
and two-dimensional systems, both for periodic lattices
and for systems with defects and random
perturbations. Our numerical results confirms the low complexity of the
ACP formulation for computing the full dynamical matrix and hence the
phonon spectrum.
%

The rest of the paper is organized as follows. Section~\ref{sec:prelim} introduces the basic formulation of KSDFT and DFPT. Section~\ref{sec:polar} describes the ACP formulation. Numerical results are presented in section~\ref{sec:numer}, followed by conclusion and discussion in section~\ref{sec:conclusion}.

\section{Preliminary}\label{sec:prelim}

For completeness we first provide a brief introduction to Kohn-Sham density
functional theory (KSDFT), and density functional perturbation theory (DFPT) in the
context of phonon calculations. To simplify our discussion, we neglect
the spin degeneracy, temperature dependence, as well as the usage of
nonlocal pseudopotential. We assume all orbitals $\{\psi_i(\vr)\}$ are
real. The spatial dimension $d=3$ is assumed in the treatment of e.g.
Coulomb interaction unless otherwise specified. We remark that such
simplified treatment does not reduce the core difficulty of the problem.

\subsection{Kohn-Sham density functional theory}
Consider a system consisting of $N_A$ nuclei and $N_e$ electrons. In
the Born-Oppenheimer approximation, for each set of nuclear
positions $\{\vR_{I}\}_{I=1}^{N_A}$, the electrons are relaxed to their ground
state. The ground state total energy is denoted by
$E_{\text{tot}}(\{\vR_{I}\}_{I=1}^{N_A})$, and can be computed
in Kohn-Sham density functional
theory~\cite{HohenbergKohn1964,KohnSham1965} according to the
minimization of the following 
Kohn-Sham energy
functional
\begin{equation}\label{eqn:KSfunc}
  \begin{split}
     E_{\KS}(\{\psi_i\};\{\vR_{I}\}) =&\frac{1}{2} \sum_{i=1}^{N_{e}} \int \abs{\nabla
    \psi_i(\vr)}^2 \ud \vr + \int V_{\text{ion}}(\vr;\{\vR_{I}\}) \rho(\vr) \ud
    \vr\\
    &
+\frac{1}{2} \iint v_{c}(\vr,\vr')\rho(\vr) \rho(\vr') \ud \vr \ud \vr' 
+ E_{\xc}[\rho] +  E_{\mathrm{II}}(\{\vR_{I}\}).
  \end{split}
\end{equation}
Here the minimization is with respect to the Kohn-Sham orbitals
$\{\psi_i\}_{i=1}^{N_e}$ satisfying the orthonormality condition
\[
\int \psi_i^*(\vr) \psi_j ( \vr) \ud \vr = \delta_{ij},\quad i,j=1,\ldots,N_e.
\]
In Eq.~\eqref{eqn:KSfunc},
$\rho(\vr) = \sum_{i=1}^{N_e}\left| \psi_i(\vr) \right|^2$
defines the electron density.   
In the discussion below we will omit the range of indices $I,i$
unless otherwise specified. 
In Eq.~\eqref{eqn:KSfunc}, $v_c(\vr,\vr') = \frac{1}{\left| \vr -
\vr' \right|}$ defines the kernel for Coulomb interaction in  $\mathbb{R}^3$. $V_{\text{ion}}$ is a local potential characterizing the electron-ion interaction in all-electron calculations, and is independent of the electronic states $\{ \psi_i\}$. More specifically, $V_{\text{ion}}$ is the summation of local potentials from each atom $I$
\begin{align}
	V_{\text{ion}}(\vr;\{\vR_{I}\}) = \sum_{I} V_I (\vr - \vR_I ).
\end{align}
In a pseudopotential approximation, it is also possible to rewrite $V_I (\vr - \vR_I )$ as
\begin{equation}
  V_{I}(\vr-\vR_{I}) = \int
  v_{c}(\vr,\vr')m_{I}(\vr'-\vR_{I}) \ud \vr',
  \label{eqn:pseudocharge}
\end{equation}
where $m_{I}$ is a localized function in the real space and is called
a pseudocharge~\cite{Martin2004,PaskSterne2005}. The normalization condition for each pseudocharge is $\int m_I(\vr) \ud \vr = -Z_I$, 
and $Z_I$ is the atomic charge for the $I$-th atom. The total pseudocharge is defined as
$m(\vr) = \sum_{I} m_{I}(\vr-\vR_{I}).$
We assume the system is charge neutral, i.e.
\[
\int m(\vr) \ud \vr = -\sum_I Z_I = -N_e.
\]
$E_{\xc}$ is the exchange-correlation energy, and here we assume
semi-local functionals such as local
density approximation (LDA)~\cite{CeperleyAlder1980,PerdewZunger1981}
and generalized gradient approximation (GGA)
functionals~\cite{Becke1988,LeeYangParr1988,PerdewBurkeErnzerhof1996} are used.
The last term in
Eq.~\eqref{eqn:KSfunc} is the ion-ion Coulomb interaction energy. For
isolated clusters in 3D, 
\begin{equation}
    E_{\mathrm{II}}(\{\vR_{I}\}) = \frac{1}{2}\sum_{I\ne J}
    \frac{Z_{I}Z_{J}}{\abs{\vR_{I}-\vR_{J}}}.
    \label{eqn:EII}
\end{equation}
We also note that for extended systems, modeled as infinite periodic structures,
both the local-pseudopotential and ion-ion terms require
special treatment in order to avoid divergence
due to the long-range $1/r$ nature of the Coulomb
interaction~\cite{Martin2004}.

The Euler-Lagrange equation associated with the Kohn-Sham energy
functional gives rise to the Kohn-Sham equations as
\begin{align}
	& H[\rho] \psi_i = \left( -\frac{1}{2} \Delta + \mcV[\rho] \right) \psi_i = \varepsilon_i \psi_i , \label{eqn:Hamil}\\
	& \int \psi_i^*(\vr) \psi_j ( \vr) \ud \vr = \delta_{ij}, \quad \rho(\vr) = \sum_{i=1}^{N_e}\left| \psi_i(\vr) \right|^2.
\end{align}
Here the eigenvalues $\{ \varepsilon_i \}$ are ordered non-decreasingly. $\psi_1,\ldots,\psi_{N_e}$ are called the occupied orbitals, while $\psi_{N_e + 1},\ldots$ are called the unoccupied orbitals.
$\psi_{N_{e}}$ is often referred to as the highest occupied
molecular orbital (HOMO), and $\psi_{N_{e}+1}$ the lowest
unoccupied molecular orbital (LUMO). The difference of the corresponding eigenvalues
$\varepsilon_{g} = \varepsilon_{N_{e}+1} - \varepsilon_{N_{e}}$ defines the HOMO-LUMO gap. Here we are interested in insulating and
semiconducting systems with positive energy gap $\varepsilon_{g}$.

For a given electron density $\rho$, the effective potential 
$\mcV[\rho]$ is
\begin{equation}
\begin{split}
\mcV [\rho](\vr) &= V_{\text{ion}}(\vr;\{\vR_{I}\}) + \int v_c(\vr,\vr') \rho(\vr') \ud \vr' + V_{\xc}[\rho](\vr) \\
&= \int v_c(\vr,\vr')(\rho(\vr')+m(\vr'))\ud \vr' + V_{\xc}[\rho](\vr). 
\end{split}
    \label{eqn:effV}
\end{equation}
Here $V_{\xc}[\rho](\vr) = \frac{\delta E_{\xc}}{\delta \rho(\vr)}$ is
the exchange-correlation potential, which is 
the functional derivative of the exchange-correlation energy with
respect to the electron density. The Kohn-Sham Hamiltonian depends
nonlinearly on the electron density $\rho$, and the electron density
should be solved self-consistently.
When the Kohn-Sham energy functional $E_{\KS}$ achieves its
minimum, the self-consistency of the electron density is simultaneously
achieved.  Then the total energy can be equivalently computed
as~\cite{Martin2004}
\begin{equation}
	\begin{split}
    E_{\text{tot}} = &\sum_{i=1}^{N_{e}}\varepsilon_{i} - \frac{1}{2}\iint 
    (\rho(\vr) - m(\vr))v_c(\vr,\vr')
    (\rho(\vr') + m(\vr')) \ud \vr \ud \vr'\\
    &-\int V_{\xc}[\rho](\vr) \rho(\vr) \ud \vr + E_{\xc}[\rho] +
    E_{\mathrm{II}}(\{\vR_{I}\}).
    \end{split}
    \label{eqn:Etotal}
\end{equation}
Here $E_{\text{band}} = \sum_{i=1}^{N_{e}}\varepsilon_{i}$
is referred to as the band energy. Using the Hellmann-Feynman
theorem~\cite{Martin2004}, the atomic force can be computed as
\begin{equation}
    \begin{split}
    \vF_{I} =& -\frac{\partial E_{\text{tot}}(\{\vR_{I}\})}{\partial
    \vR_{I}} = -\int \frac{\partial
    V_{\text{ion}}}{\partial \vR_{I}}(\vr;\{\vR_{I}\}) \rho(\vr) \ud
    \vr - \frac{\partial E_{\mathrm{II}}(\{\vR_{I}\})}{\partial \vR_{I}}\\
    =&-\int \frac{\partial
    V_{I}}{\partial \vR_{I}}(\vr-\vR_{I}) \rho(\vr) \ud
    \vr - \frac{\partial E_{\mathrm{II}}(\{\vR_{I}\})}{\partial
    \vR_{I}}.
    \end{split}
    \label{eqn:force}
\end{equation} 
The atomic force allows the performance of structural relaxation of the
atomic configuration, by minimizing the total energy $E_{\text{tot}}$
with respect to the atomic positions $\{\vR_{I}\}$. 
When the atoms are at their equilibrium positions, all atomic forces should
be $0$.

\subsection{Density functional perturbation theory}

In density functional perturbation theory (DFPT), we assume that the
self-consistent ground state electron density $\rho$ has been computed.
In this paper, we focus on phonon calculations using DFPT.
Assume the system deviates from its equilibrium position $\{\vR_{I}\}$
by some small magnitude, then the changes of the total energy is
dominated by the Hessian matrix with respect to the atomic positions. The dynamical matrix $D$ consists of 
$d\times d$ blocks in the form
\[
D_{I,J} = \frac{1}{\sqrt{M_{I}M_{J}}}\frac{\partial^2
E_{\text{tot}}(\{\vR_{I}\})}{\partial \vR_{I} \partial \vR_{J}},
\]
where $M_{I}$ is the mass of the $I$-th nuclei. The dimension of the
dynamical matrix is $d\times N_{A}$. The equilibrium atomic
configuration is at least at the local minimum of the 
total energy, and all the
eigenvalues of $D$ are real and non-negative. Hence the eigen-decomposition of $D$ is
\[
D u_{k} = \omega_{k}^2 u_{k},
\]
where $u_{k}$ is called the $k$-th phonon mode, and $\omega_{k}$ is called the
$k$-th phonon frequency.
The phonon spectrum is defined as the distribution of the eigenvalues $\{\omega_k\}$ i.e.
\begin{equation}
\varrho_D(\omega)=\frac{1}{d N_A} \sum_{k} \delta(\omega-\omega_k).
\label{eqn:phononspec}
\end{equation}
Here $\delta$ is the Dirac-$\delta$ distribution. $\varrho_D$ is also referred to as the density of states of $D$~\cite{Martin2004,LinSaadYang2016}.




In order to compute the Hessian matrix, we obtain from
Eq.~\eqref{eqn:force} that
\begin{equation}
    \begin{split}
        &\frac{\partial^2 E_{\text{tot}}(\{\vR_{I}\})}{\partial \bvec{R_{I}}
        \partial \bvec{R_{J}}} \\
         =& \int \frac{\partial V_I}{\partial
        \vR_I}(\vr-\vR_{I}) \frac{\delta \rho(\vr)}{\delta \vR_J}
        \ud\vr + \delta_{I,J}\int \rho(\vr) \frac{\partial^2 V_I}{\partial
        \vR_I^2}(\vr-\vR_{I}) \ud \vr 
        + \frac{\partial^2 E_{\mathrm{II}}(\{\vR_{I}\})}{\partial \vR_I
        \partial \vR_J}.\\
    \end{split}
    \label{eqn:SecDeriv}
\end{equation}

%
%
%

In Eq.~\eqref{eqn:SecDeriv}, the second term can be readily computed with numerical integration, and the third term involves only ion-ion interaction that is independent of the electronic states. Hence the first term is the most challenging one
due to the response of the electron density with
respect to the perturbation of atomic positions. Applying the chain
rule, we have

\begin{equation}
    \int \frac{\partial V_I}{\partial
        \vR_I}(\vr-\vR_{I}) \frac{\delta \rho(\vr)}{\delta \vR_J}
        \ud\vr = 
    \iint \frac{\partial V_I}{\partial
        \vR_I}(\vr-\vR_{I}) \frac{\delta
        \rho(\vr)}{\delta V_{\text{ion}}(\vr')}  \frac{\partial
        V_J}{\partial \vR_J}(\vr'-\vR_{J}) \ud \vr \ud \vr'.
    \label{eqn:chiterm}
\end{equation}
Here the Fr\'{e}chet derivative $\displaystyle{\chi(\vr,\vr') = \frac{\delta \rho(\vr)}{\delta V_{\text{ion}}(\vr')}}$ is referred to as
the reducible polarizability operator~\cite{OnidaReiningRubio2002}, which characterizes the
\textit{self-consistent} linear response of the electron density at
$\vr$ with respect to an external perturbation of $V_{\text{ion}}$ at
$\vr'$. However, since the Kohn-Sham equations~\eqref{eqn:Hamil} are a
set of nonlinear equations with respect to $\rho$, the self-consistent
response is still difficult to compute. Instead, the computation of
$\chi$ must be obtained through a simpler quantity 
\[
\chi_{0}(\vr,\vr') = \frac{\delta
\rho(\vr)}{\delta \mcV(\vr')},
\]
which is called the irreducible polarizability operator (a.k.a.
independent particle polarizability operator)~\cite{OnidaReiningRubio2002}. 

For simplicity in the discussion below, we will not distinguish the
continuous and discretized representation of various quantities. In
the case when a discretized representation is needed, we assume that the computational domain is uniformly discretized into a number of grid points $\{\vr_\alpha\}_{\alpha=1}^{N_g}$. We may refer to quantities such as $U(\vr)\equiv [u_1(\vr),\ldots,u_{N_e}(\vr)]$ as a matrix of dimension $N_g\times N_e$. In particular, all indices in the subscript are interpreted as the column indices of the matrix, and row indices are given in the parenthesis if necessary. For example, $u_i$ refers to the $i$-th column of $U$, and $u_i(\vr_\alpha)$ refers to the matrix element with row index $\vr_\alpha$. We denote by $M_{ij}$ a vector with a stacked column index $ij$, which refers to the $(i+(j-1)N_1)$-th column of the matrix $M$. Here the index $i$ ranges from $1$ to $N_1$, and $j$ from $1$ to $N_2$, respectively.
We also employ the following linear algebra
notation. 
The matrix-vector product $Av$ should be interpreted as
\[
(Av)(\vr) = \int A(\vr,\vr') v(\vr') \ud \vr'.
\]
Similarly the matrix-matrix product $AB$ should be interpreted as 
\[
(AB)(\vr,\vr') = \int A(\vr,\vr'') B(\vr'',\vr') \ud \vr''.
\]
The Hadamard product of two vectors $u\odot v$ should be interpreted
as
\[
(u\odot v)(\vr) = u(\vr) v(\vr).
\]


Using the chain rule and the definition of $\mcV$ in
Eq.~\eqref{eqn:effV}, we have
\begin{align}
    \chi = \frac{\delta \rho}{\delta \mcV} \left(I + \frac{\delta
    \mcV}{\delta \rho} \frac{\delta \rho}{\delta V_{\text{ion}}}\right) =
    \chi_0 ( I + v_{\text{hxc}} \chi), \label{eqn:chitmp}
\end{align}
and thus
\begin{align}
	\chi = \left(I - \chi_0 v_{\text{hxc}}\right)^{-1} \chi_0 :=
    \varepsilon^{-1}\chi_0. \label{eqn:chi}
\end{align}
Here $I$ is the identity operator, and
\[
v_{\text{hxc}} = \frac{\delta
    \mcV}{\delta \rho} =
v_{c} + \frac{\delta V_{\xc}}{\delta \rho}:= 
v_{c} + f_{\xc}.
\]
Here $f_{\xc}$ is called the exchange-correlation kernel, and is a diagonal matrix in the LDA and GGA formulation of the exchange-correlation functionals.
Eq.~\eqref{eqn:chi} also defines the operator $\varepsilon=I - \chi_0 v_{\text{hxc}}$, which is
called the dielectric operator. Using linear algebra notation, Eq.~\eqref{eqn:chiterm} requires the computation of $g_{I,a}^T\chi g_{J,b}$. 
Here $g_{I,a} = \frac{\partial V_I}{\partial {\vR_I,a}}(\vr-\vR_{I})$ is the derivative of local pseudopotential $V_I$ with respect to atomic position $R_I$ along the $a$-th direction ($a=1,\ldots,d$). Hence for phonon calculations, we need to compute $\chi g$ for $d\times N_A$ vectors in the set of $\{g_{I,a}\}$. In the discussion below, we may simply use $g$ to reflect any vector from the set  $\{g_{I,a}\}$ unless otherwise specified.

In order to compute $u=\chi g$, we apply both sides of Eq.~\eqref{eqn:chitmp} to $g$ and obtain
\begin{equation}
u = \chi_0 g + \chi_0  v_{\text{hxc}} u.
\label{eqn:dyson}
\end{equation}
Eq.~\eqref{eqn:dyson} is a fixed point problem for $u$, and is often referred to as the Dyson equation in physics literature. The simplest way to solve Eq.~\eqref{eqn:dyson} is to use a fixed point iteration
\begin{equation}
u^{k+1} = \chi_0 g + \chi_0  v_{\text{hxc}} u^{k},
\label{eqn:neumann}
\end{equation}
where $u^{k}$ is the approximate solution to $u$ at the $k$-th iteration.
The fixed point iteration~\eqref{eqn:neumann} corresponds to the Neumann expansion of the matrix inverse in Eq.~\eqref{eqn:chi}, and hence only converges when the spectral radius of $\chi_0  v_{\text{hxc}}$ is small enough (usually smaller than $1$). In order to improve the convergence behavior, more advanced numerical schemes such as Anderson mixing~\cite{Anderson1965} can be used, similar to the situation in the self-consistent field iteration for KSDFT calculations. We refer readers to~\cite{LinYang2013} for more details on solving fixed point problems in the context of electronic structure calculations.

In order to solve the Dyson equation~\eqref{eqn:dyson}, we need to apply $\chi_0$ to vectors of the form $g$ or $v_{\text{hxc}} u$. 
For systems with a finite band gap at the zero temperature, $\chi_0(\vr,\vr')$ can be computed using the Adler-Wiser formula \cite{Adler1962,Wiser1963}
\begin{align}
	\chi_0(\vr,\vr') = 2\sum_{i=1}^{N_e} \sum_{j = N_e + 1}^{\infty} \frac{\psi_i(\vr) \psi_j(\vr)  \psi_i(\vr') \psi_j(\vr' )}{\varepsilon_i - \varepsilon_j}, \label{eqn:AdlerWiser}
\end{align}
where $(\varepsilon_i, \psi_i), i=1,2,\ldots$, are the eigenpairs
in Eq.~\eqref{eqn:Hamil}. Note that $\chi_0$ is a Hermitian operator, and is negative semidefinite since $\varepsilon_i < \varepsilon_j$. In linear algebra notation, $\chi_0$ can be written as
\[
\chi_0 = 2 \sum_{i=1}^{N_e} \sum_{j = N_e + 1}^{\infty} \frac{1}{\varepsilon_i - \varepsilon_j}(\psi_i\odot \psi_j)  (\psi_i\odot \psi_j)^T.
\]

However, the Adler-Wiser formula in its original form~\eqref{eqn:AdlerWiser} can be very expensive, since it requires in principle \textit{all} the unoccupied orbitals 
$\{\psi_j\}_{j=N_e+1}^{\infty}$. These unoccupied orbitals are in general not available in KSDFT calculations, and are costly to compute and even to store in memory. However, from the Adler-Wiser formula, we can compute the multiplication of $\chi_0$ with any given vector $g$, without computing non-occupied orbitals, by solving the so called Steinheimer equation~\cite{GonzeLee1997}. This strategy has also been employed by several recent works related to the polarizability operator~\cite{UmariStenuitBaroni2010,GiustinoCohenLouie2010,NguyenPhamRoccaEtAl2012} in the context of many body perturbation theory. 
Introduce the projection operator to the unoccupied space 
\[
Q=I-\sum_{i=1}^{N_e} \psi_i\psi_i^T,
\]
then we can compute $\chi_0 g$ as
\begin{equation}
\begin{split}
\chi_0 g &= 2 \sum_{i=1}^{N_e} \sum_{j = N_e + 1}^{\infty} \frac{1}{\varepsilon_i - \varepsilon_j}(\psi_i\odot \psi_j)  (\psi_j\odot\psi_i)^T g\\
&=2\sum_{i=1}^{N_e} \psi_i\odot \left[\sum_{j = N_e + 1}^{\infty} \psi_j (\varepsilon_i - \varepsilon_j)^{-1} \psi_j^T (\psi_i \odot g)\right]\\
&=2\sum_{i=1}^{N_e} \psi_i\odot \left[Q (\varepsilon_i - H)^{-1} Q (\psi_i \odot g)\right].
\end{split}
\label{eqn:Sternheimer}
\end{equation}
Here we have used that $(\varepsilon_j,\psi_j)$ is an eigenpair of the Hamiltonian operator $H$.
Eq.~\eqref{eqn:Sternheimer} provides a practical numerical scheme for
evaluating $\chi_0 g$. The matrix inverse in Eq.~\eqref{eqn:Sternheimer}
can be avoided by solving the Steinheimer
equation~\cite{OnidaReiningRubio2002}
\[
Q(\varepsilon_i - H)Q \zeta_i =Q(\psi_i \odot g),
\]
using standard direct or iterative linear solvers. The choice of the solver can depend on practical matters such as the discretization scheme, and the availability of preconditioners. In practice for planewave discretization, we find that the use of the minimal residual method (MINRES)~\cite{PaigeSaunders1975} gives the best numerical performance.
Alg.~\ref{alg:AdlerWiser} summarizes the algorithm for computing $\chi_0 g$, without explicitly computing the unoccupied orbitals. 


\begin{algorithm}[ht]
	\small
	\DontPrintSemicolon
	\caption{Compute $\chi_0 g$.}
	\label{alg:AdlerWiser}
	
	\KwIn{
    \begin{minipage}[t]{4in}  
    Vector $g$, occupied orbitals $\{\psi_i\}_{i=1}^{N_e}$ and eigenvalues $\{\varepsilon_i \}_{i=1}^{N_e}$.
    \end{minipage}	}
	
	\KwOut{\begin{minipage}[t]{4in}  $u=\chi_0 g$. \end{minipage}
	} 
	\begin{enumerate}[leftmargin=*]
		\item Initialize $u \gets 0$. \\
		\item \textbf{For} $i=1,\ldots,N_e$ 
		\begin{enumerate}
			\item Solve the Steinheimer equation $Q(\varepsilon_i - H)Q \zeta_i =Q(\psi_i \odot g)$;
			\item $u \, \leftarrow u +  2\psi_i \odot \zeta_i$;
		\end{enumerate}
		\textbf{end for}
	\end{enumerate}
\end{algorithm}




The use of Alg.~\ref{alg:AdlerWiser} together with a proper method for
solving the Dyson equation~\eqref{eqn:dyson} gives rise to the basic
formulation for phonon calculations in DFPT. This method will be simply
referred to as DFPT for the discussion below. Assuming Dyson equations always converge in a constant number of iterations that is independent of the system size $N_e$, then main cost of DFPT is associated with the application of  $\chi_0$ to $\Or(N_e)$ vectors. Each application requires solving $N_e$ equations, and hence there are  $\Or(N_e^2)$ equations to solve. The computational cost of applying the projection operator $Q$ to a vector is $\Or(N_e^2)$, and hence the overall complexity is $\Or(N_e^4)$~\cite{BaroniGironcoliDalEtAl2001}.  Compared to the ``frozen phonon'' approach discussed in section~\ref{sec:intro}, DFPT provides an accurate description of the linear response properties of materials. 

\section{Adaptively compressed polarizability operator}\label{sec:polar}


In this section, we develop a new method for reducing the computational complexity of DFPT from $\Or(N_e^4)$ to $\Or(N_e^3)$. The reduction of the computational complexity is achieved by means of reducing the $\Or(N_e^2)$ equations in DFPT to $\Or(N_e)$ equations with systematic control of the accuracy. In particular, our method does not employ the ``nearsightedness'' property of electrons for insulating systems with substantial band gaps as in linear scaling methods~\cite{Kohn1996}. Hence our method can be applied to insulators as well as semiconductors with small band gaps. In section~\ref{sec:prelim}, we have reduced the problem of computing the dynamic matrix to the computation of $\chi g_{I,a}$, where $\{g_{I,a}\}$ is a set of fixed vectors given by the derivative of the local pseudopotential with respect to the atomic positions.  Let us stack the indicies $I,a$ into a single index $j$, and denote by 
\begin{align}
	G := \left[g_{1},\ldots,g_j,\ldots,g_{d\times N_A}\right]
\end{align}
the matrix collecting all these vectors. More generally, $G$ can be any fixed matrix with $\Or(N_e)$ columns as required in different applications of DFPT. Then our method consists of two main steps: 1) Find a compressed representation of $\chi_0$, which allows the computation of $\chi_0 G$ by solving only $\Or(N_e)$ linear equations. 2) Update the compressed representation of $\chi_0$, which allows the accurate computation of $\chi G$ without significant increase of the computational cost. In particular, step 2) requires the compression strategy of $\chi_0$ to be \textit{adaptive} to the solution of the Dyson equation~\eqref{eqn:dyson}. Hence we refer our representation of $\chi_0$ as the adaptively compressed polarizability operator (ACP). The steps 1) and 2) of the ACP formulation are given in section~\ref{sec:compresschi0} and~\ref{sec:compresschi}, respectively.




\subsection{Compression of $\chi_0$}\label{sec:compresschi0}

Consider first the computation of $\chi_0 G$ as required in the initial step in Eq.~\eqref{eqn:dyson}. In general, the singular values of $\chi_0$ decay slowly, and a forcefully applied low rank decomposition of $\chi_0$ such as those based on the singular value decomposition (SVD) will lead to inaccurate results. Nonetheless, it is possible to find a compressed representation of $\chi_0$ when we only need to evaluate $\chi_0 G$ for a fixed matrix $G$. 

According to  Alg.~\ref{alg:AdlerWiser}, computing $\chi_0 G$ involves solving the following $\Or(N_e^2)$ Steinheimer equations 
\begin{equation}
	Q(\varepsilon_i - H)Q \zeta_{ij} = Q(\psi_i \odot g_j), \quad i = 1, \ldots,N_e,\quad j=1,\ldots,d\times N_A.
	\label{eqn:eqnU}
\end{equation}
As $N_e$ becomes large, asymptotically there can be many more equations
to solve than the dimension of the matrix $N_g\sim \Or(N_e)$, and hence
it should be possible to compress the redundant information in the right
hand side vectors. In fact this observation has been used in various
contexts in computational chemistry for compressing the Hadamard product
of occupied and unoccupied orbitals, which is called ``density fitting''
(DF) or ``resolution of identity'' (RI) techniques to compress
$\Or(N_e^2)$ vectors into $\Or(N_e)$ vectors with a relatively small
pre-constant~\cite{Weigend2002,RenRinkeBlumEtAl2012}. 

It should be noted that  density fitting techniques alone do not reduce the number of equations to solve. The reason is that the Eqs.~\eqref{eqn:eqnU} have the dependence on the shift $\varepsilon_i$ on the left hand side. Hence even if the number of right hand side vectors is reduced to $\Or(N_e)$, multiplied with the $N_e$ shifts, we still have $\Or(N_e^2)$ equations to solve! Therefore, in order to reduce the complexity for computing $\chi_0 G$, we must disentangle the right hand side vectors and the shifts. Note that all $\{\varepsilon_i\}$ are eigenvalues corresponding to occupied orbitals, and are typically contained in a relatively small interval (in the order of eV), at least in the pseudopotential framework. 

More specifically, consider the following parameterized equation
\begin{equation}
	Q(\varepsilon - H)Q \zeta = \xi,
	\label{eqn:eqnparam}
\end{equation}
where $\xi$ is any vector in the range of $Q$. Since $\varepsilon \in \mc{I}\equiv [\varepsilon_1, \varepsilon_{N_e}]$ , we can systematically obtain the solution to the parameterized equation by evaluating on a few sampled points in $\mc{I}$. In this work, we choose the Chebyshev nodes $\{\wt{\varepsilon}_c\}_{c=1}^{N_c}$, which are obtained by a linear map the Chebyshev nodes in the reference interval $[-1,1]$ to $\mc{I}$, i.e.
\[
\wt{\varepsilon}_c = \frac{\varepsilon_1+\varepsilon_{N_e}}{2} + \frac{\varepsilon_{1}-\varepsilon_{N_e}}{2} \cos \theta_c, \quad \theta_c = \frac{\pi(c-\frac12)}{N_c}, \quad c=1,\ldots,N_c.
\]
Typically it is sufficient to choose the number of Chebyshev nodes $N_c$ to be $10\sim 40$. Denote by $\wt{\zeta}_c$ the solution to Eq.~\eqref{eqn:eqnparam} corresponding to $\varepsilon=\wt{\varepsilon}_c,\quad c=1,\ldots,N_c$, then any solution $\zeta$ with $\varepsilon\in \mc{I}$ can be obtained by a Lagrange interpolation procedure as


\begin{equation} \label{eqn:LagPoly}
\zeta = \sum_{c=1}^{N_c} 
\wt{\zeta}_c \prod_{c'\ne c} \frac{\varepsilon - \wt{\varepsilon}_{c'}}{\wt{\varepsilon}_c - \wt{\varepsilon}_{c'}}.
\end{equation}

Using Chebyshev interpolation~\eqref{eqn:LagPoly}, we need to solve Eq.~\eqref{eqn:eqnU} with $\varepsilon_i$ replaced by $\wt{\varepsilon}_c$. At first sight, the number of equations does not decrease but actually increase by a factor of $N_c$ compared to the original formulation~\eqref{eqn:eqnU}. However, Chebyshev interpolation disentangles the index $i$ that appears both in the shift and in the right hand side. Since $N_c$ is a constant that is independent of the system size, if we can find a compressed representation of the right hand side vectors using $\Or(N_e)$ vectors, we reduce the overall number of equations to solve to $\Or(N_e)$.

Let us denote by $M$ the collection of right hand side vectors in Eq.~\eqref{eqn:eqnU} without the $Q$ factor. More specifically, 
\[
M_{ij} = \psi_i\odot g_j,\quad \mbox{or}\quad M_{ij}(\vr) = \psi_i(\vr)g_j(\vr).
\]
Here we have used $ij$ as a stacked column index for the matrix $M$.
The dimension of $M$ is $N_g\times \Or(N_e^2)$. Typically, the
computational complexity for the compression for such a dense matrix $M$
with approximate rank $\Or(N_e)$ is $\Or(N_e^4)$, even with the help of
the recently developed randomized algorithms
(see~\cite{HalkoMartinssonTropp2011} for a review). Nonetheless, note
that
for a fixed row index $\vr_\alpha$, the row vector given by $\{M_{ij}(\vr_\alpha)\}_{i,j=1}^{i=N_e,j=d\times N_A}$ is the Kronecker product between the row vector given by 
$\{g_j(\vr_\alpha)\}_{j=1}^{d\times N_A}$ and that given by
$\{\psi_i(\vr_\alpha)\}_{i=1}^{N_e}$. As will be seen below, this structure allows the computational complexity of the compression of $M$ to be reduced to $\Or(N_e^3)$.

To this end, we seek for the following interpolative decomposition (ID) type of compression~\cite{ChengGimbutasMartinssonEtAl2005} for the matrix $M$, i.e.
\begin{equation}\label{eqn:MID}
M_{ij}(\vr) \approx \sum_{\mu=1}^{N_\mu} \xi_\mu(\vr) M_{ij}(\vr_\mu) \equiv  \sum_{\mu=1}^{N_\mu} \xi_\mu(\vr) \psi_i(\vr_\mu) g_j(\vr_\mu).
\end{equation}
Here $\{\vr_\mu\}_{\mu=1}^{N_\mu}$ denotes a collection of selected row
indices. Numerical results indicate that it is sufficient to choose
$N_\mu\sim \Or(N_e)$ with a relatively small pre-constant.
Interpolative decomposition has been developed in the context of
compressing the four center integrals in quantum chemistry calculations 
with $\Or(N_e^3)$ computational complexity~\cite{LuYing2015}.
Here we extend this method to compress $M$.

The first step of the randomized density fitting method is to employ a
fast pre-processing procedure, such as a subsampled random Fourier transform
(SRFT)~\cite{WoolfeLibertyRokhlinEtAl2008}, to transform the matrix $M$
into a smaller matrix $\wt{M}$ of dimension $N_g\times \Or(N_e)$. 
The second step is to use a pivoted QR decomposition to obtain an interpolative decomposition of $\wt{M}$, which gives the interpolation vectors 
$\{\xi_\mu\}$ as well as the selected row indices $\{\vr_\mu\}$ as needed in Eq.~\eqref{eqn:MID}.   
We summarize the procedure for compressing $M$ in Alg.~\ref{alg:compressM}, and refer readers to ~\cite{LuYing2015} for more details of the algorithm.

\begin{algorithm}[ht]
	\small
	\DontPrintSemicolon
	\caption{Interpolative decomposition for $M$ using a randomized
  density fitting method~\cite{LuYing2015}.}
	\label{alg:compressM}
	
	\KwIn{    
    \begin{minipage}[t]{4in}
    Matrix $M$. Threshold tolerance $\epsilon$. 
    \end{minipage}\\
     
	}
	
	\KwOut{\begin{minipage}[t]{4in}  
    Selected row indices $\{\vr_\mu\}$, and interpolation vectors $\{\xi_\mu\}$.
    \end{minipage}
	}
	
	\begin{enumerate}[leftmargin=*]
		\item Subsampled random Fourier transform of $M$:
		\begin{enumerate}
			\item Compute for $\nu = 1,\ldots, N_e\times d\times N_A$ the discrete Fourier transform
			\[
			\hat{M}_\nu (\vr) = \sum_{I=1}^{N_e\times d\times N_A} e^{ -2\pi \I I \nu / (N_e\times d\times N_A) } \eta_I M_I(\vr),
			\]
			where $\eta_I$ is a random complex number with unit modulus for each $I$.
			\item Choose a submatrix $\wt{M}$ of matrix $\hat{M}$ by randomly choosing $rN_e$ columns. In practice, $r = 8$ and $r = 16$ are used in our implementation for one-dimensional and two-dimensional numerical examples, respectively.
		\end{enumerate}
		\item Compute the pivoted QR decomposition of the $rN_e \times N_g$ matrix $\wt{M}^T$ : $\wt{M}^T \wt{\Pi} = \wt{Q}\wt{R}$, where the absolute values of the diagonal entries of $\wt{R}$ are ordered non-increasingly.
		\item Determine the number of selected columns $N_\mu$, such that $|\wt{R}_{N_\mu+1,N_\mu+1} | < \epsilon |\wt{R}_{1,1}|  \le | \wt{R}_{N_\mu,N_\mu} | $. Form $\{\vr_\mu\}$, $\mu = 1,\ldots,N_\mu$ such that the $\vr_\mu$-column of $\wt{M}^T$ corresponds to the $\mu$-th column of $\wt{M}^T \wt{\Pi}$.
	    \item Denote by $\wt{R}_{1:N_\mu,1:N_\mu}$ the submatrix of $\wt{R}$ consisting of its first $N_\mu \times N_\mu$ entries, and $\wt{R}_{1:N_\mu,:}$ the submatrix consisting of the first $N_\mu$ rows of $\wt{R}$. Compute
	    \[
      \Xi^{T} = \wt{R}_{1:N_\mu,1:N_\mu}^{-1} \wt{R}_{1:N_\mu,:} \wt{\Pi}^{-1}.
	    \]
	    Then the $\mu$-th column of the $N_g \times N_\mu$ matrix $\Xi$
      gives the interpolation vector $\xi_\mu$.
	\end{enumerate}
   
\end{algorithm}

\begin{remark}
In Alg.~\ref{alg:compressM}, step 1.(a), it is possible to avoid the explicit construction of the matrix $M$. Instead of performing SRFT on the entire matrix $M$, we could apply SRFT only to the matrix $G$, and select $r$ columns as a matrix $\wt{G}$.  Then for a fixed row index $\vr_\alpha$, the Kronecker product between the rows of subsampled matrix $\wt{G}$, $\{\wt{g}_i(\vr_\alpha)\}_{i=1}^{r}$, and $\{\psi_i(\vr_\alpha)\}_{i=1}^{N_e}$ gives one row for $\wt{M}$. In practice we find that this heuristic procedure also works well for compressing the matrix $M$ in phonon calculations.
\end{remark}

Once the compressed representation~\eqref{eqn:MID} is obtained, we solve the following set of modified Steinheimer equations
\begin{equation}\label{eqn:eqncompressU}
Q(\wt{\varepsilon}_c - H)Q \wt{\zeta}_{c\mu} = Q \xi_\mu,\quad c=1,\ldots,N_c,\quad \mu=1,\ldots,N_\mu.
\end{equation}
Here $c\mu$ is the stacked column index for $\wt{\zeta}$. The
 number of equations is hence reduced to
$N_c N_{\mu}\sim \Or(N_e)$. Using Eq.~\eqref{eqn:LagPoly}, we construct the  quantity $W = \left[ W_1,\ldots,W_{N_\mu} \right]$. Each column of $W$ is defined by

\begin{equation}
W_{\mu} = 2\sum_{i=1}^{N_e}  \psi_i \odot \left(\sum_{c=1}^{N_c} 
\wt{\zeta}_c \prod_{c'\ne c} \frac{\varepsilon_i - \wt{\varepsilon}_{c'}}{\wt{\varepsilon}_c - \wt{\varepsilon}_{c'}}\right) \psi_i(\vr_\mu). 
\label{eqn:eqnW}
\end{equation}
Combining Eq.~\eqref{eqn:eqnW} with Eq.~\eqref{eqn:Sternheimer}, we obtain directly $\chi_0 g_j $ as
\begin{equation}
\chi_0 g_j \approx \sum_{\mu=1}^{N_\mu} W_\mu g_j(\vr_\mu).
\label{eqn:chi0gcompress}
\end{equation}
It should be noted that in Eq.~\eqref{eqn:chi0gcompress}, we have avoided the explicit reconstruction of the solution vectors $\zeta_{ij}$ as in Eqs.~\eqref{eqn:eqnU}, of which the computational cost is again $\Or(N_e^4)$.  

Formally, Eq.~\eqref{eqn:chi0gcompress} can further be simplified by
defining a matrix $\Pi$ with $N_\mu$ columns, which consists of selected
columns of a permutation matrix, \ie~ 
$\Pi = \wt{\Pi}_{:,1:N_\mu}$ 
as the first $N_\mu$ columns of the permutation matrix obtained
from pivoted QR decomposition. More specifically,
$\Pi_{\mu}=e_{\vr_\mu}$ and $e_{\vr_\mu}$ is a unit vector with only one
nonzero entry at $\vr_\mu$ such that $e_{\vr_\mu}^T g_j=g_j(\vr_\mu)$.
Then 
\begin{equation}
\chi_0 g_j \approx W \Pi^T g_j := \wt{\chi}_0 g_j.
\label{eqn:chi0gcompress2}
\end{equation}
Note that the compressed polarizability operator $\wt{\chi}_0=W \Pi^T$ is  
 formally independent of the right hand side vector $\{g_j\}$, and the rank of $\wt{\chi}_0$ is only $N_\mu$, while the singular values of $\chi_0$ have a much slower decay rate. This is because $\wt{\chi}_0$ only agrees with $\chi_0$ when applied to vectors $g_j$. In other words, the difference between $\wt{\chi}_0$ and $\chi_0$ is not controlled in the space orthogonal to that spanned by $G$.  Alg.~\ref{alg:chi_0} summarizes the algorithm for computing the compressed polarizability operator $\wt{\chi}_0$.

\begin{algorithm}[ht]
	\small
	\DontPrintSemicolon
	\caption{Computing compressed polarizability operator $\wt{\chi}_0$.}
	\label{alg:chi_0}
	
	\KwIn{
    \begin{minipage}[t]{4in}  
			Vectors $\{g_j\}$. 
    Hamiltonian matrix $H$. 
    
    Eigenpairs corresponding to occupied orbitals 
    $\{\psi_i,\varepsilon_i\}$
    \end{minipage}\\
	}
	
	\KwOut{\begin{minipage}[t]{4in}  
    $\wt{\chi}_0 = W\Pi^T$.
    \end{minipage}
	} 
	\begin{enumerate}[leftmargin=*]
		\item Use Alg.~\ref{alg:compressM} to obtain $\{\vr_\mu\}$, $\Pi$, and hence the compressed representation of $M$. 
		\item Solve compressed Eqs.~\eqref{eqn:eqncompressU}.
		\item Compute $W$ using Eq.~\eqref{eqn:eqnW}.
	\end{enumerate}
	
\end{algorithm}

	
			
	
	

The computational complexity of Alg.~\ref{alg:chi_0} can be analyzed as follows. For simplicity we neglect all possible logarithmic factors in the complexity analysis. The cost for constructing the compressed representation of $M$ is $\Or(N_e^3)$. Eqs.~\eqref{eqn:eqncompressU} require solving $N_cN_\mu\sim \Or(N_e)$ equations.  Assuming the computational cost for applying $H$ to a vector is $\Or(N_g)$, and assuming that the number of iterations using an iterative solver to solve Eqs.~\eqref{eqn:eqncompressU} is bounded by a constant, then the cost for solving all equations is dominated by the computation of $\{Q\xi_\mu\}$ which is $\Or(N_e^3)$. 
In order to construct $W$, for each $\mu$ and $i$, we can first compute the term in the parenthesis in the right hand side of Eq.~\eqref{eqn:eqnW}. Then the computational complexity for constructing $W$ is again $\Or(N_e^3)$. Therefore, the overall asymptotic computational cost for constructing the compressed polarizability operator $\wt{\chi}_0$ is $\Or(N_e^3)$.  In practice, we find that the computational cost is dominated by solving the $\Or(N_e)$ linear equations in step 2 of Alg.~\ref{alg:chi_0}.



\subsection{Compression of $\chi$}\label{sec:compresschi}

According to Eq.~\eqref{eqn:dyson}, $\chi_0 G$ is the leading order approximation to $U=\chi G$, and this approximation can be  inaccurate if $\chi_0 v_{\mathrm{hxc}}$ is not small. From the perspective of section~\ref{sec:compresschi0}, the self-consistent solution to the Dyson equation~\eqref{eqn:dyson} introduces two additional difficulties: 1) 
 we need to find compressed representation $\wt{\chi}_0$ that agrees
 with $\chi_0$ when applied to both $G$ and $v_{\mathrm{hxc}}U$; 2) $U$
 is not known \textit{a priori}.  Hence if we apply Alg.~\ref{alg:chi_0}
 directly, we may need to increase the rank of $\wt{\chi}_0$ to
 $2N_\mu$ or higher to maintain the accuracy. Below we introduce the adaptively compressed polarizability operator (ACP) method that simultaneously addresses the above two difficulties. 

We assume that $v_{\mathrm{hxc}}$ is invertible, and $v_{\mathrm{hxc}}^{-1} g$ for a vector $g$ can be computed easily. This is the case in the absence of the exchange-correlation kernel $f_{\xc}$, and $v_{\mathrm{hxc}}^{-1}g = v_c^{-1} g$ can  simply be obtained by applying the Laplacian operator to $g$. This approximation is referred to as the ``random phase approximation'' (RPA) in physics literature~\cite{OnidaReiningRubio2002}. In the presence of $f_{\xc}$ in the LDA and GGA formulation, $f_{\xc}$ is a diagonal matrix, and $v_{\mathrm{hxc}}^{-1}g = v_c^{-1} g$ can be solved using iterative methods. 

We introduce the following change of variable
\begin{equation}
U = \wt{U} - B, \quad B = v_{\mathrm{hxc}}^{-1}G,
\end{equation}
and the Dyson equation~\eqref{eqn:dyson} becomes
\begin{equation}
\wt{U} = \chi_0 v_\text{hxc} \wt{U} + B.
\label{eqn:dyson2}
\end{equation}
The advantage of using Eq.~\eqref{eqn:dyson2} over \eqref{eqn:dyson} is that formally, we only need to find $\wt{\chi}_0$ that is accurate when applied to $v_\text{hxc} \wt{U}$. In an iterative algorithm, for a given matrix $\wt{U}$, we can use Alg.~\ref{alg:chi_0} to construct $\wt{\chi}_0[\wt{U}]$ by replacing $G$ with  $v_\text{hxc} \wt{U}$, with $\wt{U}$ in the bracket to highlight the $\wt{U}$-dependence of the compression scheme, i.e.
\begin{equation}
\wt{U} = \wt{\chi}_0[\wt{U}] v_\text{hxc} \wt{U} + B.
\label{eqn:dyson3}
\end{equation}
We note that when self-consistency is reached for Eq.~\eqref{eqn:dyson3} with the self-consistent solution denoted by $\wt{U}^*$,  $\wt{\chi}_0[\wt{U}^*] v_\text{hxc} \wt{U}^*$ remains a good approximation to  $\chi_0 v_\text{hxc} \wt{U}^*$, even if $\wt{U}^*$ deviates away from the initial guess. In each step, the approximate rank of $\wt{\chi}_0[\wt{U}]$ remains to be $N_\mu$. Hence $\wt{\chi}_0[\wt{U}]$ is adaptive to the solution $\wt{U}$, and hence is called the adaptively compressed polarizability operator (ACP). This concept of adaptively constructing a low rank matrix shares similar spirit to the recently developed adaptively compressed exchange operator (ACE) for the efficient solution of Hartree-Fock-like calculations~\cite{Lin2016ACE}. 



Eq.~\eqref{eqn:dyson3} can be solved using the fixed point iteration or more advanced methods for solving fixed point problems, similar to that in Eq.~\eqref{eqn:dyson} in DFPT. However, thanks to the low rank structure of $\wt{\chi}_0$ in Eq.~\eqref{eqn:chi0gcompress2}, we can significantly accelerate the convergence. Let us denote the value of $\wt{U}$ at the $k$-th iteration as $\wt{U}^{k}$, which gives rise to the ACP 
$\wt{\chi}_0[\wt{U}^k] = W^{k} (\Pi^{k})^T$.
Eq.~\eqref{eqn:dyson2} indicates that if the magnitude of $\chi_0$ is small, then $\wt{U}^0=B$ is a good initial guess to start the iteration. 
Then we can reformulate Eq.~\eqref{eqn:dyson3} and obtain the following iteration scheme
\begin{equation}
\label{eqn:dyson4}
\wt{U}^{k+1} = \left( I - W^{k}(\Pi^{k})^T v_\text{hxc} \right)^{-1} B = B +  W^{k} \left(I-(\Pi^{k})^T v_\text{hxc} W^{k}\right)^{-1}(\Pi^{k})^T v_\text{hxc} B.
\end{equation}
The second equality in Eq.~\eqref{eqn:dyson4} uses the Sherman-Morrison-Woodbury identity for computing the inverse. The cost of the inversion is $\Or(N_e^3)$ due to the low rank structure of $\wt{\chi}_0[\wt{U}^k]$. Numerical results indicate that the iteration scheme~\eqref{eqn:dyson4} can converge much more rapidly compared to the fixed point iteration for Eq.~\eqref{eqn:dyson}. In fact often two to four iterations are sufficient to obtain results that are sufficiently accurate. 
Alg.~\ref{alg:chigACP} describes the algorithm for using ACP to compute $\chi G$.

\begin{algorithm}[ht]
  \small
  \DontPrintSemicolon
  \caption{Computing $\chi G$ with adaptively compressed polarizability operator}
  \label{alg:overall}

	\KwIn{    
    \begin{minipage}[t]{4in}  \end{minipage}\\
    Vectors $\{g_j\}$.  Stopping criterion $\delta$.
    
    Eigenpairs corresponding to occupied orbitals 
    $\{\psi_i,\varepsilon_i\}$ 
	}
   \KwOut{
   \begin{minipage}[t]{4in}  
   $U\approx \chi G$
   \end{minipage}
   } 
  \begin{enumerate}[leftmargin = *]
  	\item Compute $\wt{U}^0= B = v_\text{hxc}^{-1} G$. $k\gets 0$.
  	\item \textbf{Do}  
  	\begin{enumerate}
  		\item Use Alg.~\ref{alg:chi_0} by replacing $G$ with $v_\text{hxc} \wt{U}^k$ to obtain $W^k$ and $\Pi^k$, and obtain $\wt{\chi}_0^k = W^k (\Pi^k)^T$. \\
   		\item Update $\wt{U}^{k+1}$ according to Eq.~\eqref{eqn:dyson4}.\\
  		\item $k\leftarrow k+1$
  	\end{enumerate}
  	\textbf{until} $\|\wt{U}^{k} - \wt{U}^{k-1}\| < \delta$ or maximum number of iteration is reached.\\
  	\item Compute $U \leftarrow \wt{U}^{k} - B$.
  \end{enumerate}
  \label{alg:chigACP}
\end{algorithm}

\section{Numerical examples}\label{sec:numer}

In this section, we demonstrate the performance of ACP proposed in the
previous section, and compare it with the density functional
perturbation theory (DFPT),
and with the finite difference approach (FD) through three examples. The
first example consists of a one-dimensional (1D) reduced Hartree-Fock
model problem that can be tuned to resemble an insulating or a
semi-conducting system. The second example is a two-dimensional (2D)
model problem with a periodic triangular lattice structure. The third example is
a 2D triangular lattice with defects and random
perturbations of the atomic positions.  All results are performed on a single computational
core of a 1.4 GHz processor with 256 GB memory using MATLAB.


\subsection{One-dimensional reduced Hartree-Fock model}
The 1D reduced Hartree-Fock model was introduced by Solovej
\cite{Solovej1991}, and has been used for analyzing defects in solids
in e.g. \cite{CancesDeleurenceLewin2008,CancesDeleurenceLewin2008a}. The
simplified 1D model neglects the contribution of the
exchange-correlation term. As discussed in previous sections, the
presence of exchange-correlation functionals at LDA/GGA level   does not
lead to essential difficulties in phonon calculations.

The Hamiltonian in our 1D reduced Hartree-Fock model is given by 
\begin{equation}
	H [ \rho ] = -\frac{1}{2} \frac{d^2}{dx^2} + \int K(x,y) \left( \rho(y) + m(y) \right) \ud y.
	\label{eqn:1DHamil}
\end{equation}

Here $m(x) = \sum_{I} m_I(x - R_I)$ is the summation of pseudocharges. Each function $m_I(x)$ takes the form of a one-dimensional Gaussian 
\begin{align}
m_I(x) = - \frac{Z_I}{\sqrt{2\pi\sigma_I^2}} \exp \left( -\frac{x^2}{2\sigma_I^2}  \right),
\end{align}
where $Z_I$ is an integer representing the charge of the $I$-th nucleus. In our numerical simulation, we choose all $\sigma_I$ to be the same.

Instead of using a bare Coulomb interaction which diverges in 1D when $x$ is large, we use a Yukawa kernel as the regularized Coulomb kernel
\begin{equation}
K(x, y) = \frac{2\pi e^{-\kappa | x-y | }}{\kappa \epsilon_0}, \label{eqn:yukawa}
\end{equation}
which satisfies the equation
\begin{align}
-\frac{d^2}{dx^2}K(x,y) + \kappa^2 K(x,y) = \frac{4\pi}{\epsilon_0}\delta(x-y).
\end{align}
As $\kappa \rightarrow 0$, the Yukawa kernel approaches the bare Coulomb interaction given by the Poisson equation. The parameter $\epsilon_0$ is used so that the magnitude of the electron static contribution is comparable to that of the kinetic energy.
The ion-ion repulsion energy $E_{\mathrm{II}}$ is also computed using the Yukawa interaction $K$ in the model systems.





The parameters used in the model are chosen as follows. Atomic units are used throughout the discussion unless otherwise mentioned. For all systems tested in this subsection, the distance between each atom and its nearest neighbor is set to 2.4 a.u.. The Yukawa parameter $\kappa = 0.1$. The nuclear charge $Z_I$ is set to 1 for all atoms, and $\sigma_I$ is set to be 0.3. The Hamiltonian operator is represented in a plane wave basis set. 

By adjusting the parameter $\epsilon_0 = 1.0$ or $10$, the reduced
Hartree-Fock model can be tuned to resemble an insulator or a semiconductor, respectively. We apply ACP to both cases. We use Anderson mixing for SCF iterations, and the linearized eigenvalue problems are solved by using the locally optimal block preconditioned conjugate gradient (LOBPCG) solver~\cite{Knyazev2001}. 

\begin{figure}[!ht]
	\begin{minipage}[h]{\linewidth}
			\begin{minipage}[t]{0.5\linewidth}
				\centering
				\includegraphics[width=2.2in]{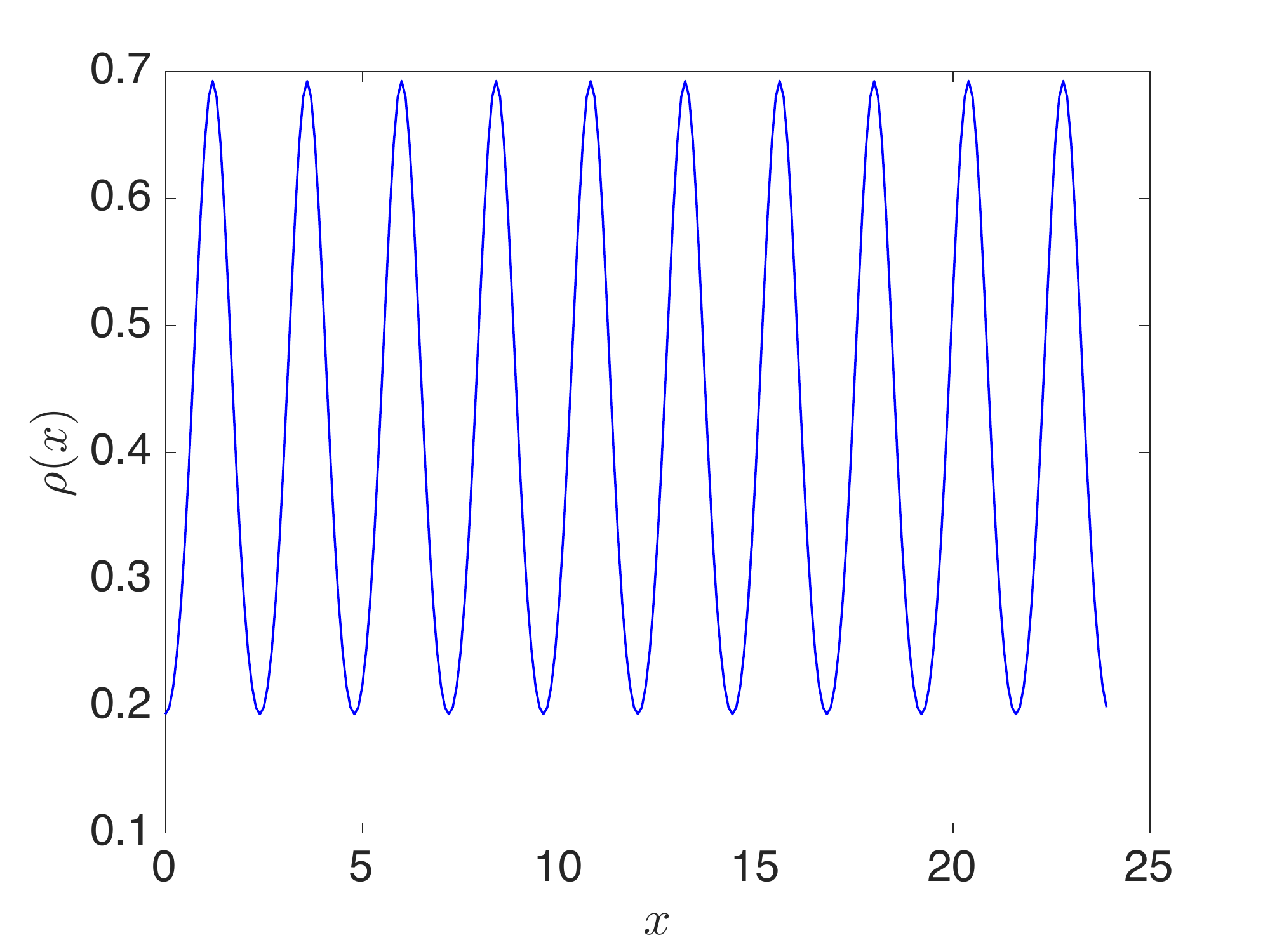}
				\subcaption{Insulator}
				\label{fig:rho60}
			\end{minipage}%
			\begin{minipage}[t]{0.5\linewidth}
				\centering
				\includegraphics[width=2.2in]{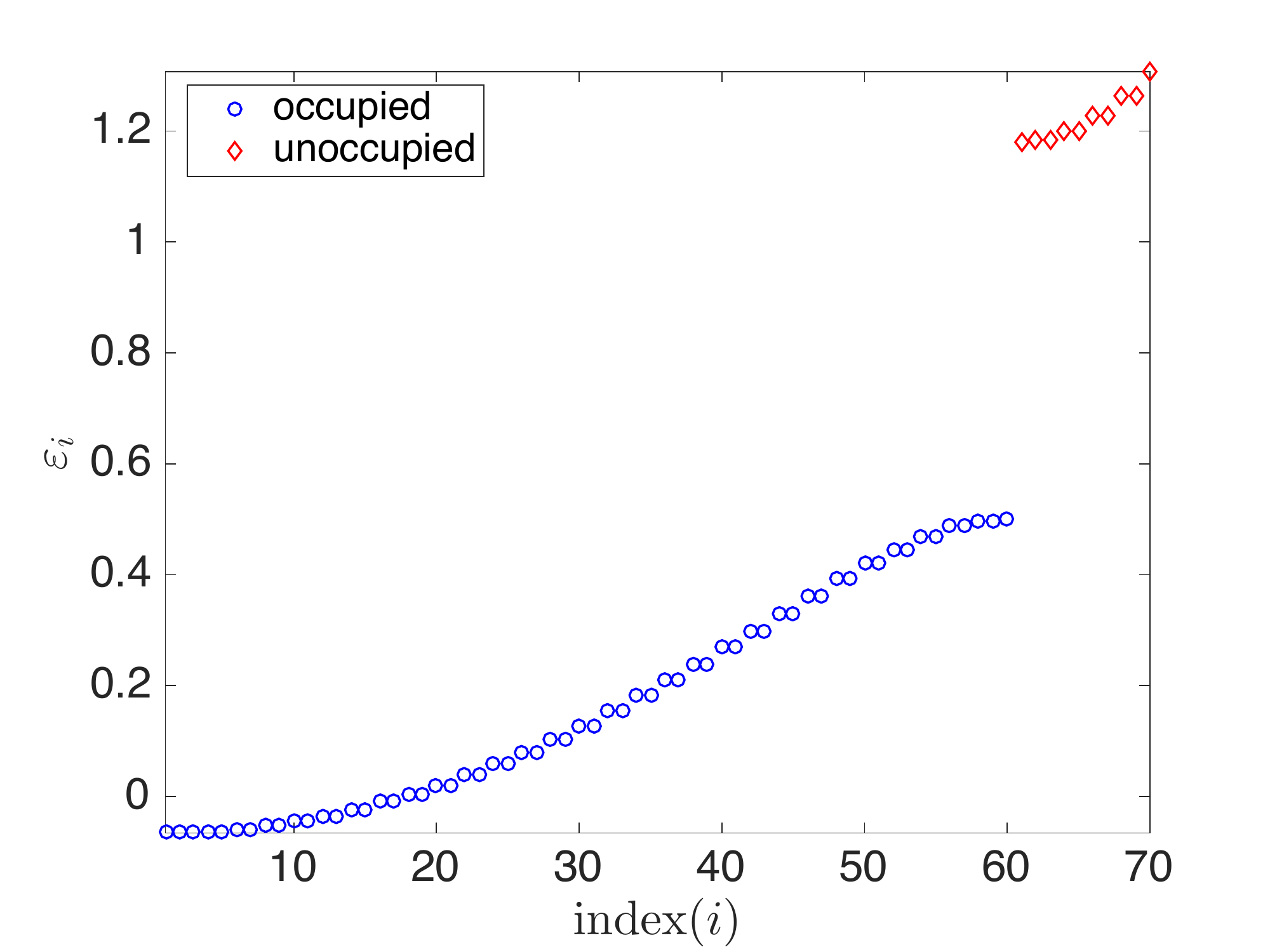}
				\subcaption{Insulator}
				\label{fig:gap60}
			\end{minipage}
	\end{minipage}
    \begin{minipage}[h]{\linewidth}
	\begin{minipage}[h]{0.5\linewidth}
		\centering
		\includegraphics[width=2.2in]{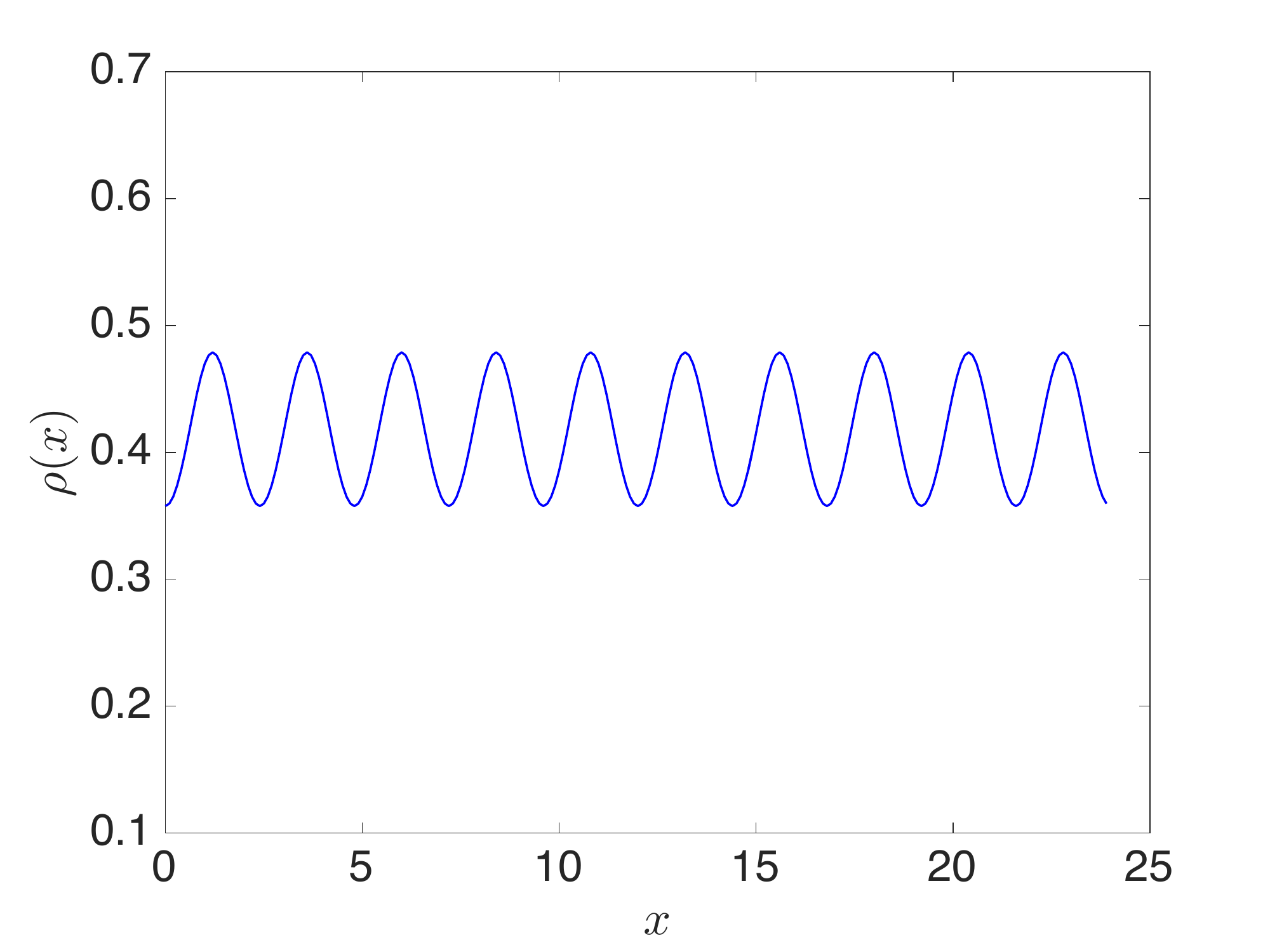}
		\subcaption{Semiconductor}
		\label{fig:rho60semi}
	\end{minipage}%
	\begin{minipage}[h]{0.5\linewidth}
		\centering
		\includegraphics[width=2.2in]{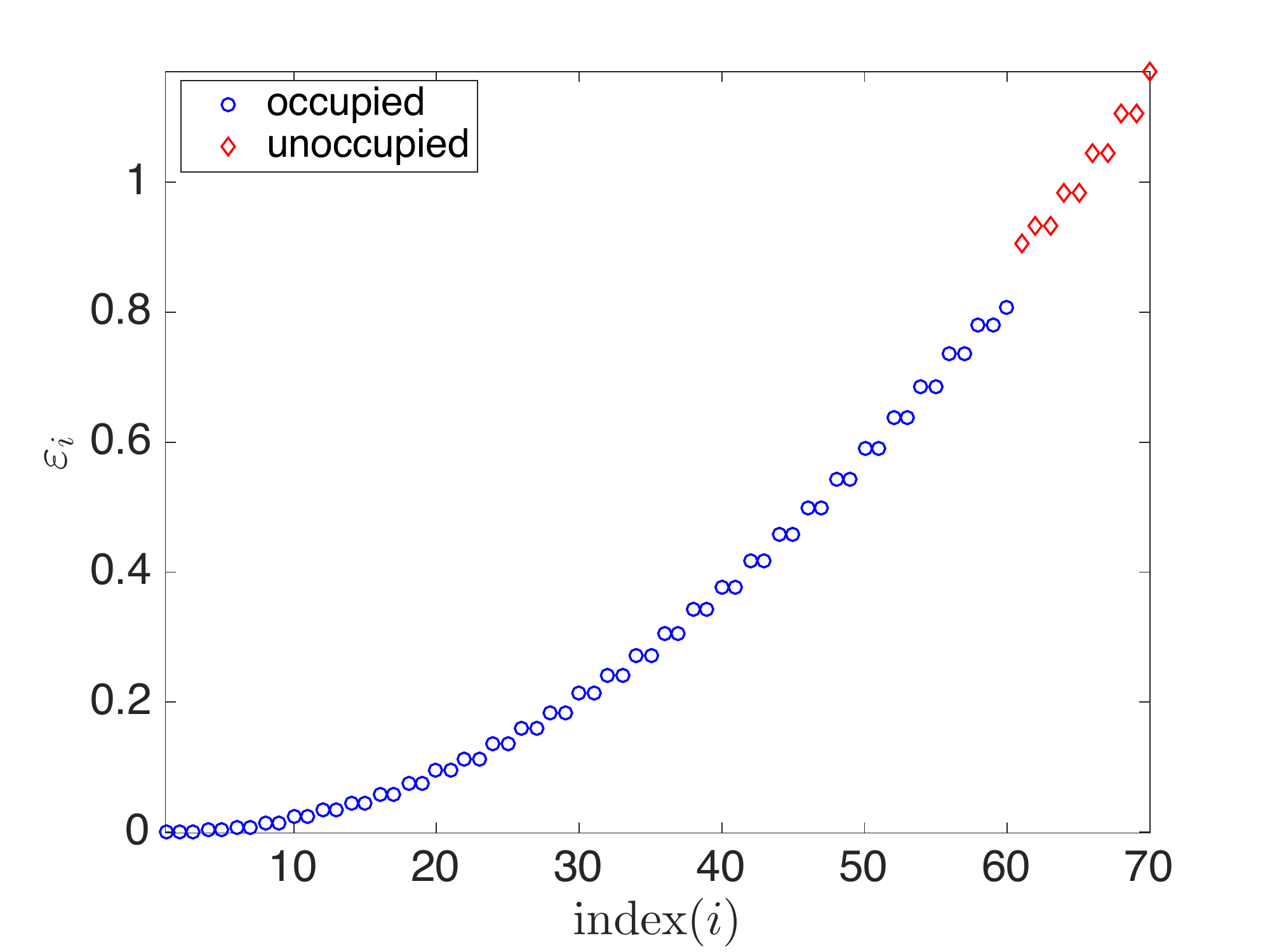}
		\subcaption{Semiconductor}
		\label{fig:gap60semi}
	\end{minipage}
    \end{minipage}
  \caption{The electron density $\rho(x)$ of the 60-atom (a) insulating,
  and (c) semiconducting systems in the left panel. The corresponding
  occupied (blue circles) and unoccupied eigenvalues (red diamonds) are
  shown in the right panel in (b), (d), respectively.}
	\label{fig:1D}
\end{figure}

For systems of size $N_A = 60$, the converged electron density $\rho$ associated with the two 1D test cases as well as the 70 smallest eigenvalues associated with the Hamiltonian defined by the converged $\rho$ are shown in Fig.~\ref{fig:1D}. For the insulator case, the  electron density fluctuates between 0.1935 and 0.6927. There is a finite HOMO-LUMO gap, $\varepsilon_g = \varepsilon_{61} - \varepsilon_{60} = 0.6763$. The electron density associated with the semiconductor case is relatively uniform in the entire domain, with the fluctuation between 0.3576 and 0.4788. The corresponding band gap is 0.1012. 
Fig.~\ref{fig:1D} is obtained by a system with 60 atoms, and we find that systems with different sizes shows similar patterns in the band structure for both insulating and semiconducting systems, respectively.

All numerical results of the ACP method below are benchmarked with
results obtained from DFPT. In order to demonstrate the effectiveness of the ACP formulation for compressing $U=\chi G$, we directly measure the relative $L^2$ error, defined as $\norm{U - U^{\text{ACP}}}_2 / \norm{U}_2$, where $ U^{\text{ACP}}$ is obtained from Alg.~\ref{alg:overall}. We also report the error of the phonon spectrum by computing the $L^\infty$ error of the phonon frequencies $\{\omega_k\}$. Due to the presence of acoustic phonon modes for which $\omega_k$ is close to 0, we report the absolute error instead of the relative error for the phonon frequencies.
In DFPT, we use
MINRES~\cite{PaigeSaunders1975} to solve the Steinheimer equations
iteratively.  The initial guess vectors for the solutions 
are obtained from previous iterations in the Dyson equation to reduce
the number of matrix-vector multiplications. The same strategy for
choosing the initial guess is implemented for the ACP
formulation as well.
Anderson mixing is used to
accelerate the convergence of Dyson equations in DFPT. In ACP we find
that the fixed point iteration~\eqref{eqn:dyson4} is sufficient for
fast convergence. 


\begin{table}[!ht]
	\begin{center}
		\begin{tabular}{c|c|c|c|c} \hline
			\diagbox{$N_c$}{$N_\mu$} & $3N_e$ & $4N_e$ & $5N_e$ & $6N_e$\\ \hline
			5 & 0.0231& 0.0018 & 0.0017 & 0.0017 \\ \hline
			10 & 0.0221 & 5.9734e-04 & 2.2727e-05 & 2.0089e-05 \\ \hline
			15 & 0.0190 & 5.9914e-04 & 1.1439e-05 & 1.5539e-05  \\ \hline
			20 & 0.0217 & 7.2434e-04 & 1.5262e-05 & 7.5518e-06 \\ \hline
		\end{tabular}
	\end{center}
	\caption{The relative $L^{2}$ error of $U=\chi G$ obtained from ACP
  compared to DFPT.
  }
	\label{tab:accuracy1}
\end{table}
In Table \ref{tab:accuracy1}, we calibrate the accuracy of our algorithm
with different number of Chebyshev points $N_{c}$ and different number of
columns $N_{\mu}$. 
We choose $N_\mu = l N_e$, where $l = 3,4,5,6.$
Table \ref{tab:accuracy1} shows that with
fixed Chebyshev interpolation points $N_c$, the numerical accuracy increases monotonically with respect to $N_\mu$, until limited by the accuracy of the Chebyshev interpolation procedure.
Similarly the increase of Chebyshev interpolation reduces the error
until being limited by the choice of $N_{\mu}$.
When both $N_{\mu}$ and $N_{c}$ are large enough, the error of
$\chi G$ can be less than $10^{-5}$.


\begin{table}[!ht]
	\begin{center}
		\begin{tabular}{c|c} \hline
			Method and parameter & $L^{\infty}$-norm error   \\ \hline
			FD $\delta = 0.01$ & 5.6779e-04\\ \hline
			ACP $N_\mu = 241$ for $\epsilon = 10^{-3}$ & 3.6436e-04 \\ \hline
			ACP $N_\mu = 359$ for $\epsilon = 10^{-5}$ & 2.7380e-06  \\ \hline
		\end{tabular}
	\end{center}
  \caption{$L^{\infty}$ error of the phonon spectrum. System is insulating with size $N_A = 60$. Chebyshev  nodes $N_c=20$. $N_\mu$ is determined such that $|\wt{R}_{N_\mu+1,N_\mu+1} | < \epsilon |\wt{R}_{1,1}|  \le | \wt{R}_{N_\mu,N_\mu}| $ in Alg.~\ref{alg:compressM}.
}
	\label{tab:accuracy2}
\end{table}

In Table \ref{tab:accuracy2}, we choose $N_\mu$ based on the entries of
$\tilde{R}$ as is shown in Alg.~\ref{alg:compressM}, and compare the results
to those obtained from the FD. In the FD approach, we set the convergence tolerance for LOBPCG to be $10^{-6}$, and the SCF tolerance to be $10^{-8}$. $\delta =0.01$ denotes the deviation of atom positions to their equilibrium ones. We remark that the varying $\delta$ from $0.01$ to $0.0001$ does not change too much in the phonon spectrum. The same parameters for SCF and LOBPCG are used to converge the ground state calculation in the ACP formulation. The absolute error of the phonon spectrum is smaller than $ 10^{-3}$.  The ACP formulation can lead to very accurate phonon
spectrum by solving a relatively small number of equations.

\begin{figure}[!ht]
	\begin{center}
		\includegraphics[width=0.5\textwidth]{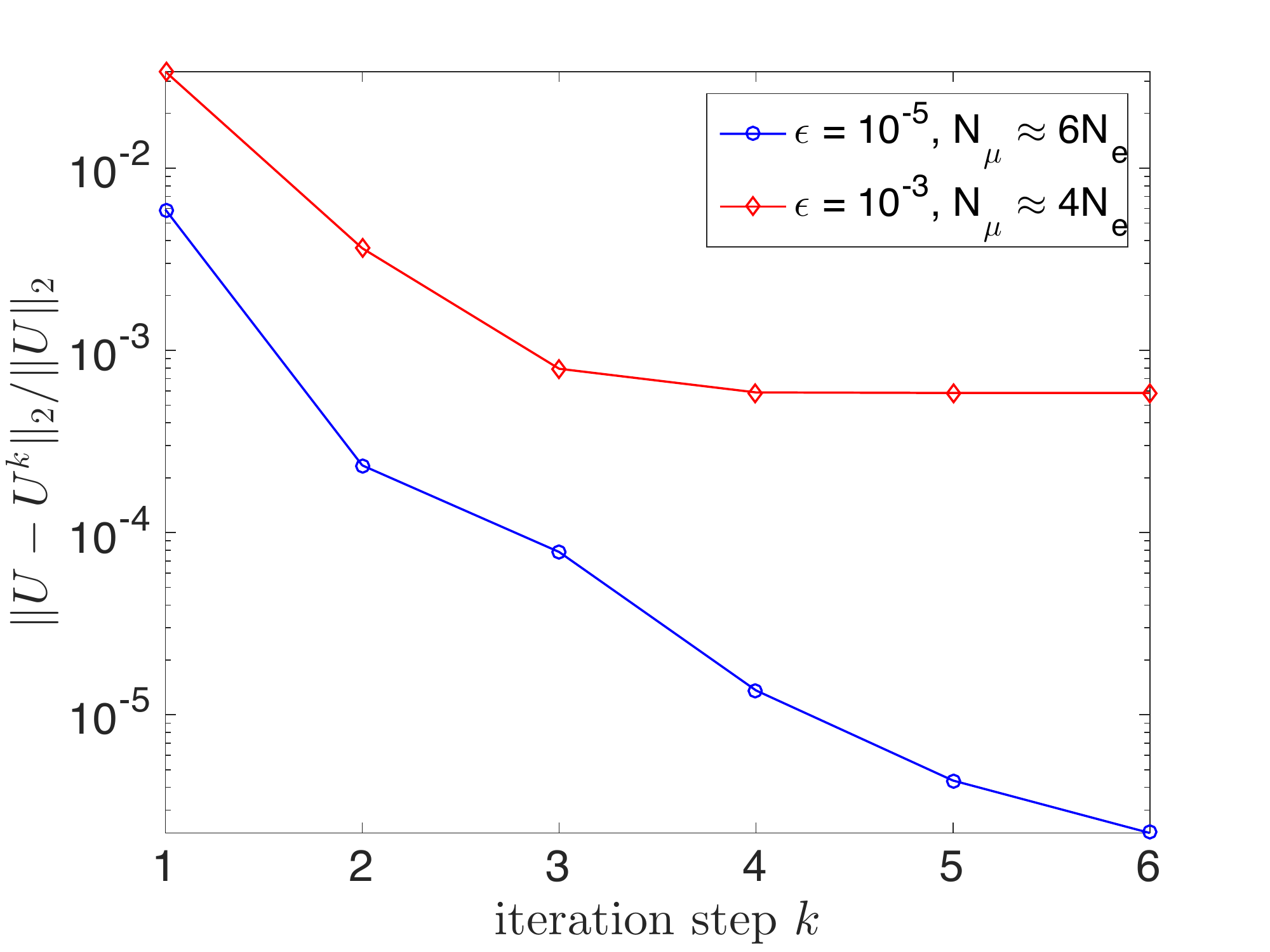}
	\end{center}
	\caption{Convergence of adaptive compression.}
	\label{fig:adapconv}
\end{figure}

In order to demonstrate the effectiveness of the adaptive compression
strategy,
the relative $L^2$ error for the approximation of $U=\chi G$ with respect to the iteration in Alg.~\ref{alg:overall} is given in Fig.~\ref{fig:adapconv}. For $\epsilon = 10^{-5}$, the error is around $10^{-5}$ after 4 iterations. For $\epsilon = 10^{-3}$, the error is reduced to $0.0008$. In this case, if we stop after the first iteration in Alg.~\ref{alg:overall}, the relative error of $\chi G$ is $0.0339$. Numerical results shows significant improvement after two to four iterations. This indicates that the self-consistent solution of the Dyson equation is crucial for the accurate computation of phonon spectrum. 

\begin{figure}[!ht]
	\begin{minipage}[t]{0.5\linewidth}
		\centering
		\includegraphics[width=2.2in]{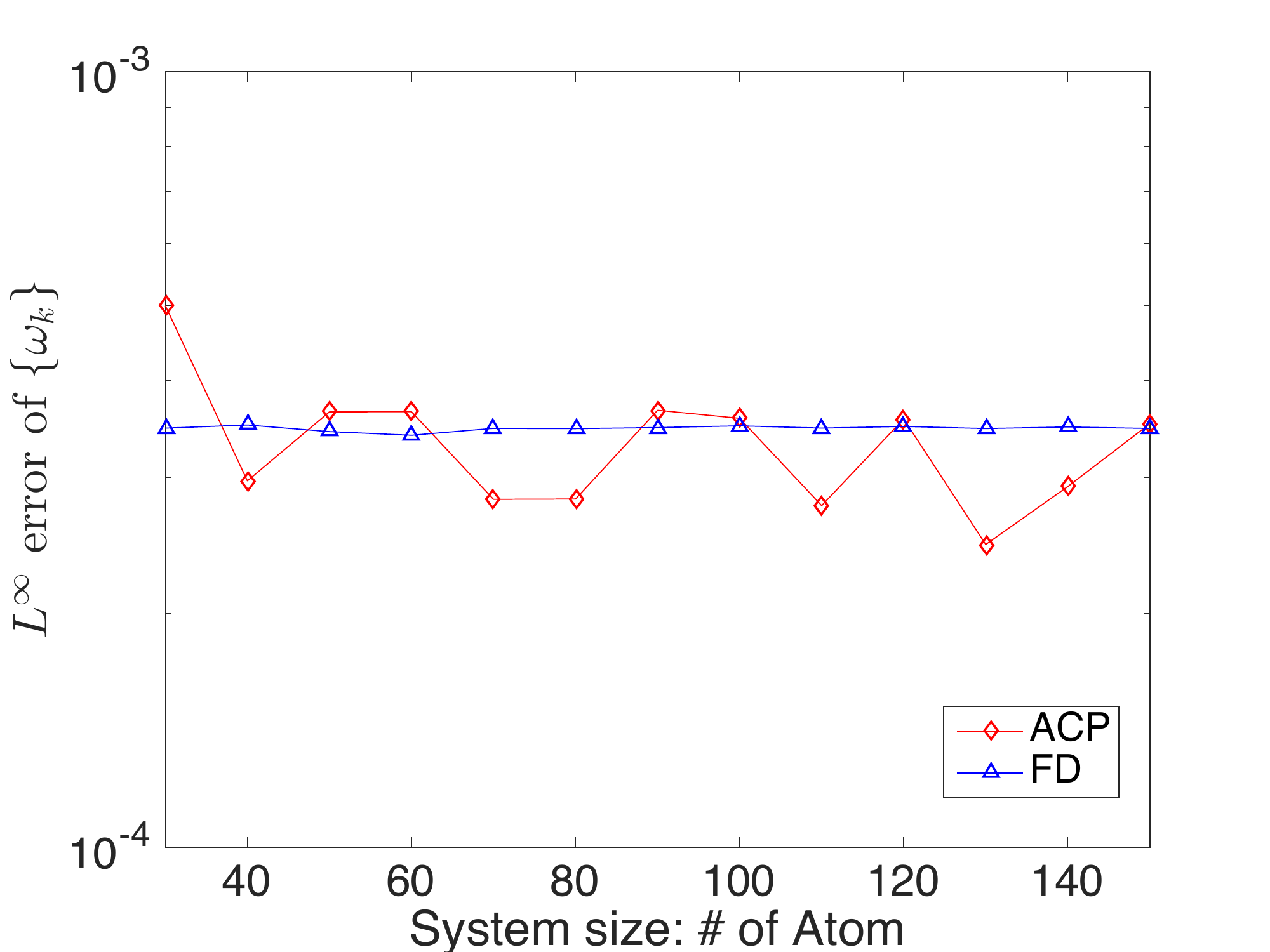}
		\subcaption{Insulator}
		\label{fig:errInsul}
	\end{minipage}%
	\begin{minipage}[t]{0.5\linewidth}
		\centering
		\includegraphics[width=2.2in]{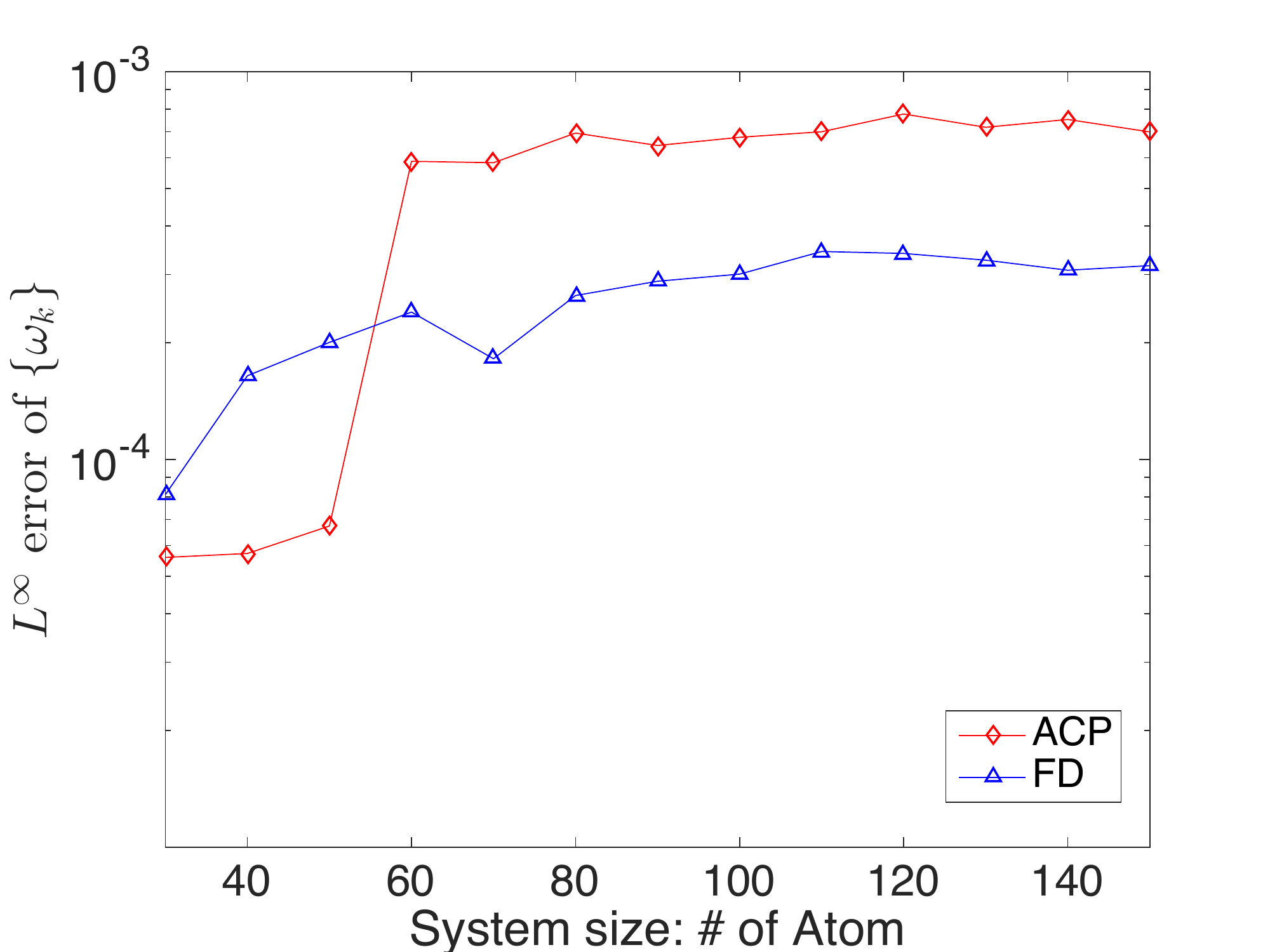}
		\subcaption{Semiconductor}
		\label{fig:errSemi}
	\end{minipage}
	\caption{$L^{\infty}$ error of the phonon frequencies $\{\omega_k\}$ obtained from ACP and FD. For ACP formulation $N_c = 20$. $N_\mu \approx 4N_e$.}
	\label{fig:err1D}
\end{figure}

\begin{figure}[!ht]
	\begin{minipage}[t]{0.5\linewidth}
		\centering
		\includegraphics[width=2.2in]{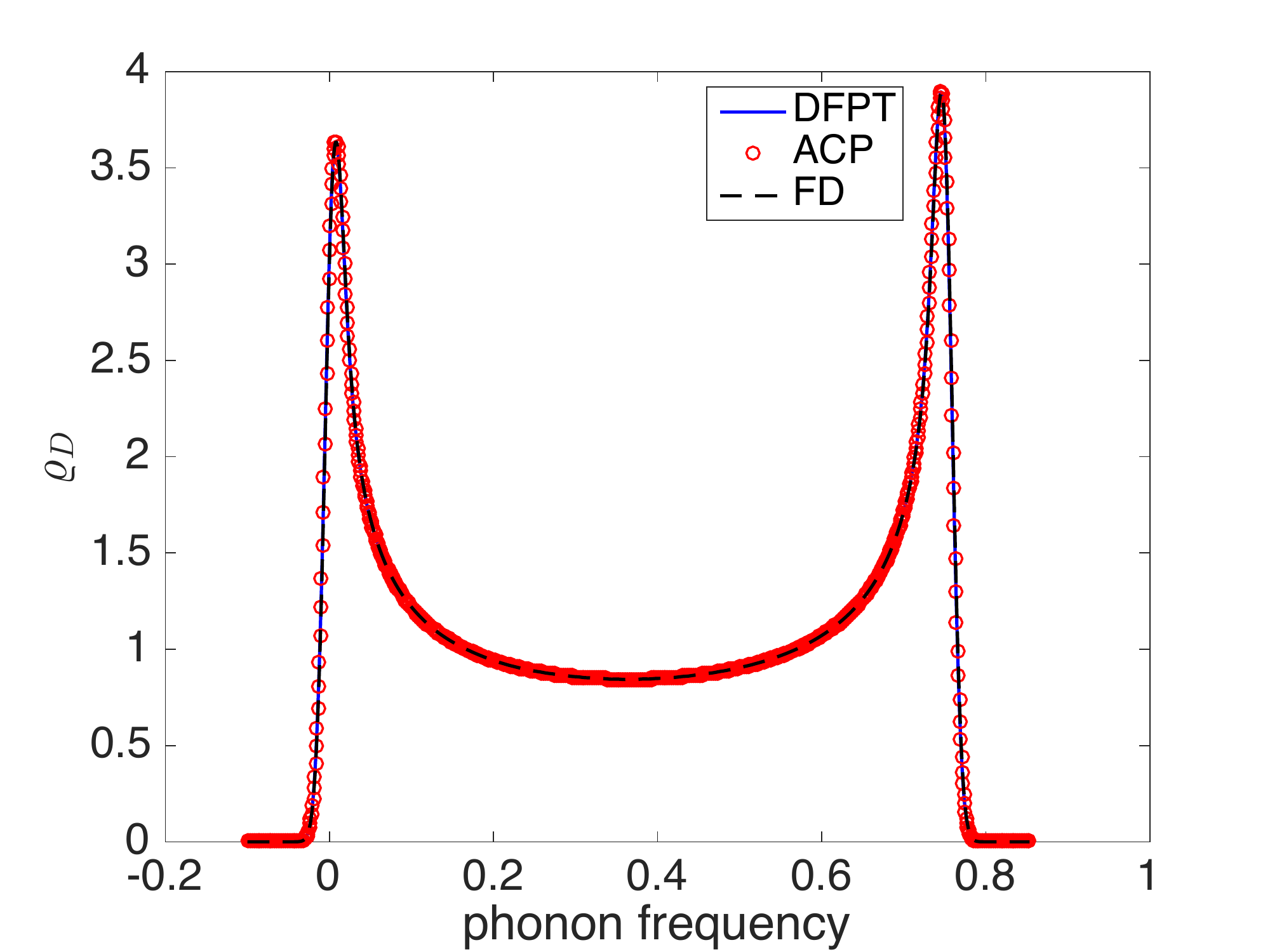}
		\subcaption{Insulator}
		\label{fig:phononspec}
	\end{minipage}%
	\begin{minipage}[t]{0.5\linewidth}
		\centering
		\includegraphics[width=2.2in]{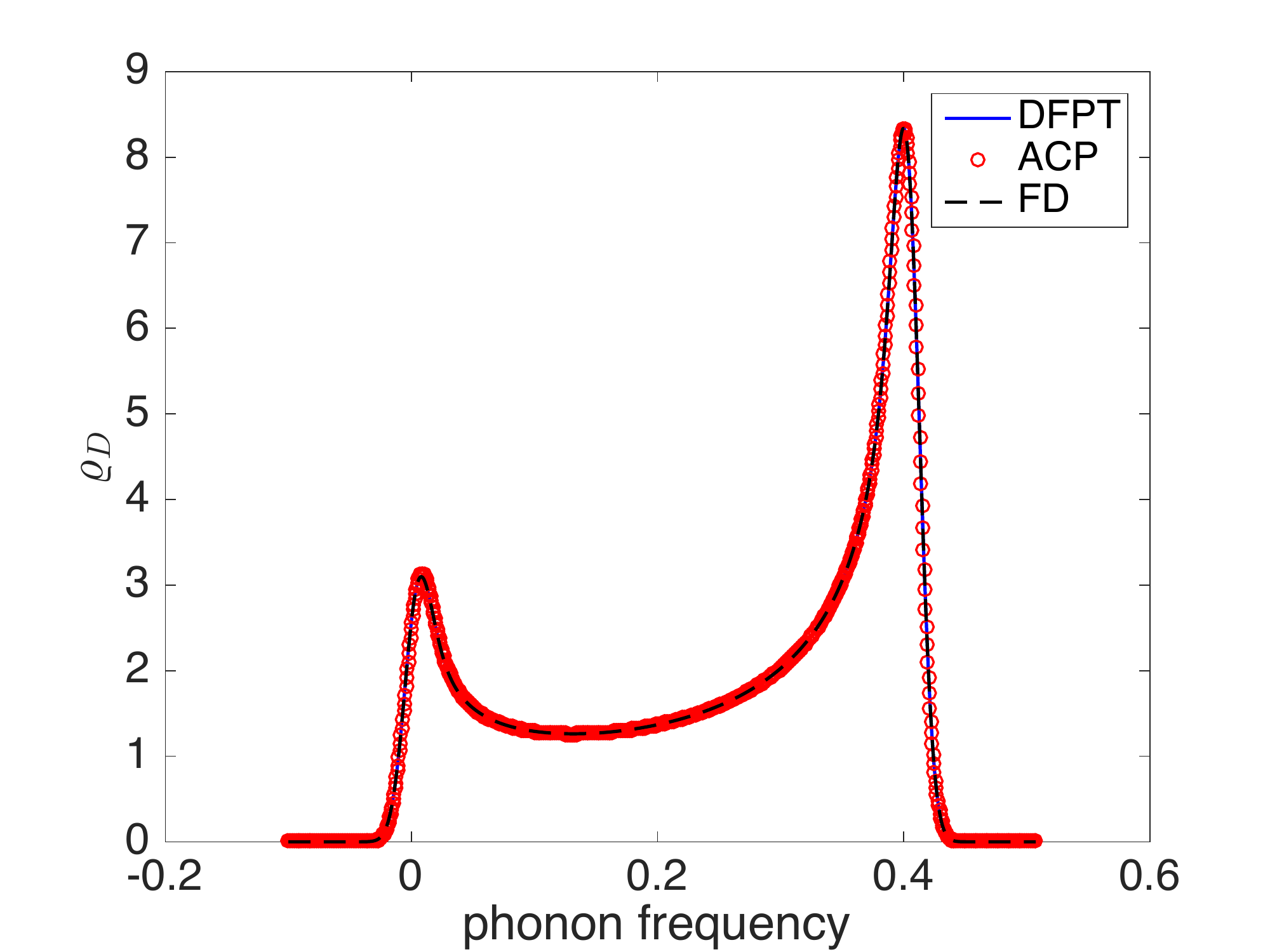}
		\subcaption{Semiconductor}
		\label{fig:phononspec_semi}
	\end{minipage}
	\caption{Phonon spectrum for the 1D systems computed using ACP, DFPT, and FD, for both (a) insulating and (b) semiconducting systems.}
	\label{fig:phononspec1D}
\end{figure}

\begin{figure}[!ht]
	\begin{minipage}[t]{0.5\linewidth}
		\centering
		\includegraphics[width=2.2in]{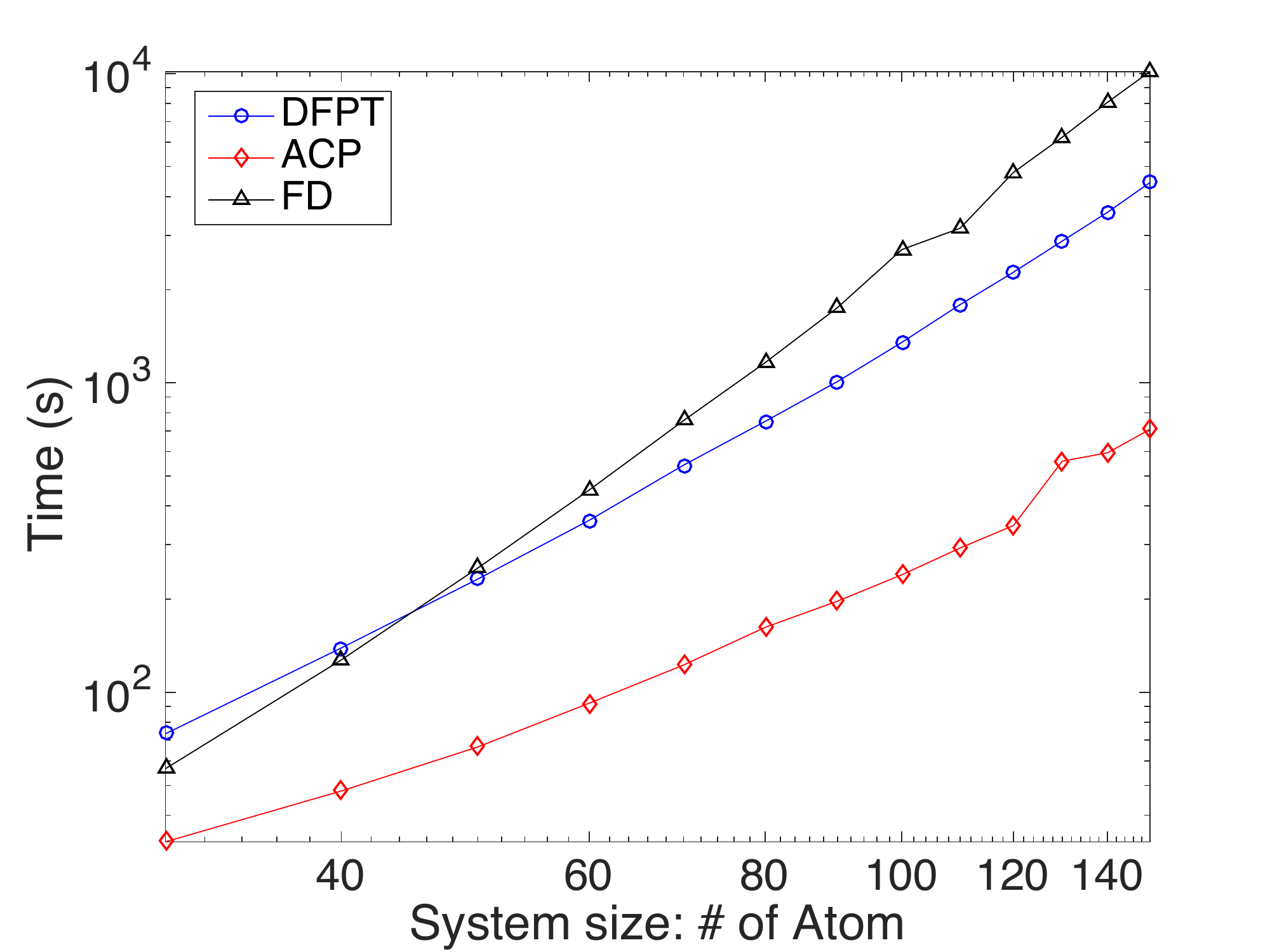}
		\subcaption{Insulator}
		\label{fig:eff_1e-3}
	\end{minipage}%
	\begin{minipage}[t]{0.5\linewidth}
		\centering
		\includegraphics[width=2.2in]{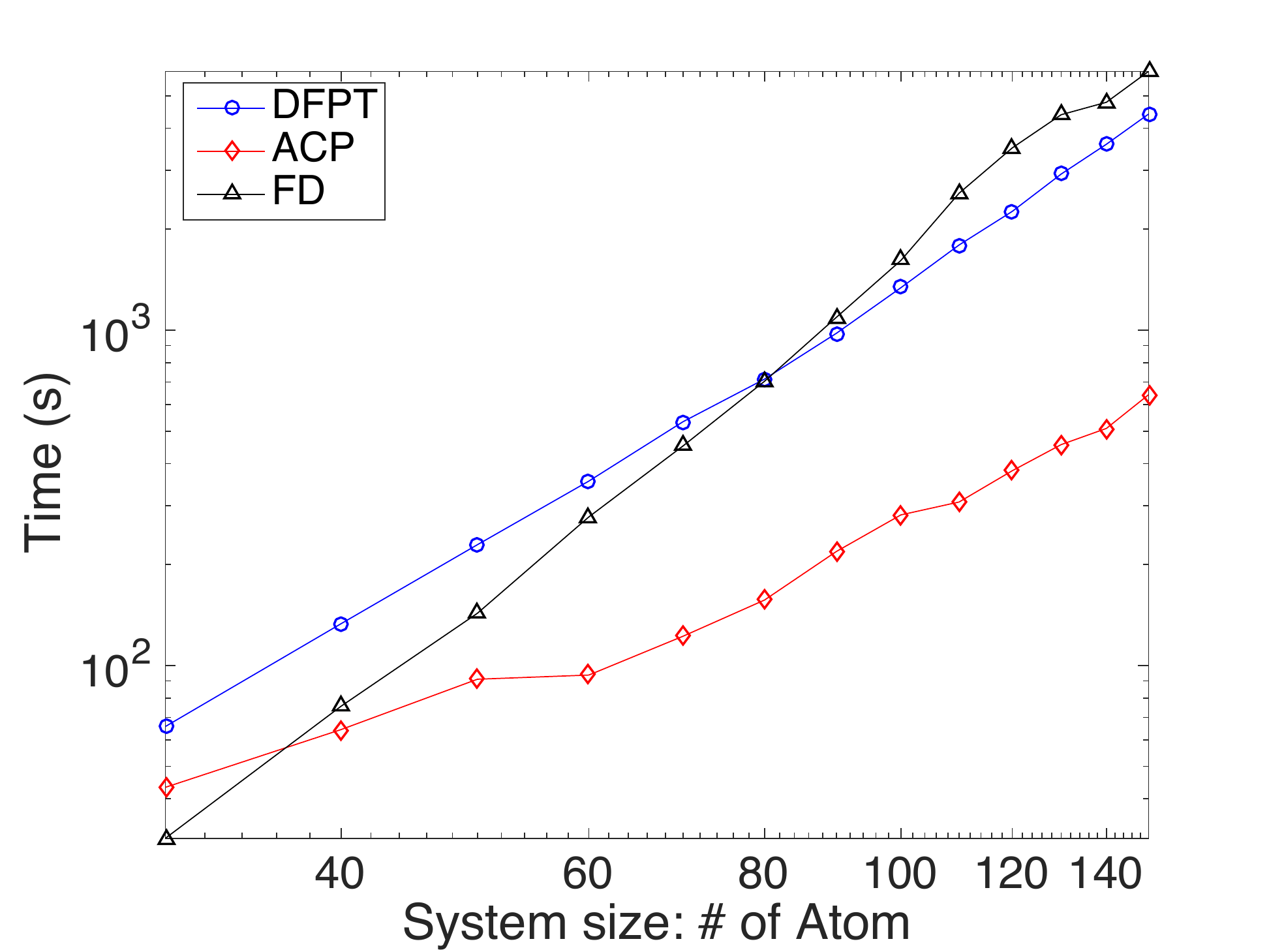}
		\subcaption{Semiconductor}
		\label{fig:eff_1e-3_semi}
	\end{minipage}
	\caption{Computational time of 1D examples. Comparison among DFPT, ACP, and FD for (a) insulating, and (b) semiconducting systems, respectively.}
	\label{fig:eff1D}
\end{figure}

We perform phonon calculations for systems of size from 30 to 150 for both insulating and semiconducting systems. In terms of accuracy, Fig.~\ref{fig:err1D} shows that as $N_\mu$ increases linearly with respect to
the system size, the accuracy of phonon spectrum ($L^{\infty}$ error)
remains to be roughly the same, which is empirically around the $\epsilon=10^{-3}$ used to determine $N_\mu$ in Alg.~\ref{alg:compressM}. 
For the computation of the phonon frequency, we find that $N_c=20$ and $N_\mu = 4N_e$ is sufficient to achieve error around $10^{-3}$, and the phonon spectrum is already indistinguishable from that obtained from
DFPT. 
Fig.~\ref{fig:phononspec1D} reports the phonon spectrum $\varrho_D$ for the systems of size $N_A = 150$. We remark that Fig.~\ref{fig:phononspec1D} plots the density of states $\varrho$ by replacing the Dirac-$\delta$ distribution in Eq.~\ref{eqn:phononspec} with regularized delta function
\[
\delta_\sigma(x) = \frac{1}{\sqrt{2\pi \sigma^2}} e^{-\frac{x^2}{2\sigma^2}}.
\]
Here the smear parameter $\sigma$ is chosen as 0.01.


To demonstrate the efficiency of the ACP algorithm, 
Fig.~\ref{fig:eff1D} compares the computational time of ACP,
DFPT, and FD, respectively.  We observe that the computational cost of DFPT matches to that of FD due to the choice of initial guess of Steinheimer equations and the Anderson mixing strategy for solving the Dyson equation. Compared to DFPT, the ACP formulation benefits both from that it solves less number of Steinheimer equations, and from that the Sherman-Morrison-Woodbury procedure~\eqref{eqn:dyson4} is more efficient than Anderson mixing for solving the Dyson equation. In fact for all systems, 
Alg.~\ref{alg:overall} converges within $4$ iterations, while the Anderson mixing in DFPT may require $20$ iterations or more for systems of all sizes.
Hence for both insulating and semiconducting systems, the ACP formulation becomes more advantageous than DFPT and FD for systems merely beyond $40$ atoms.
For the largest system with $150$ atoms, ACP is 6.28 and 6.87 times faster than DFPT for insulating and semiconducting systems, respectively.



\begin{table}[!ht]
	\begin{center}
		\begin{tabular}{c|c|c} \hline
			Method & Insulator & Semiconductor   \\ \hline
			FD  & 3.4403 & 3.3047 \\ \hline
			DFPT  & 2.8997 & 2.9459 \\ \hline
			ACP  & 2.5040 & 2.1065 \\ \hline
		\end{tabular}
	\end{center}
	\caption{Computational scaling measured from $N_A = 90$ to $N_A = 150$. }
	\label{tab:slope1D}
\end{table}

Table \ref{tab:slope1D} measures the slope of the computational cost with respect to system sizes from $N_A = 90$ to $N_A = 150$. In theory, the asymptotic computational cost of FD and DFPT should be $\Or(N_e^4)$, and the cost of ACP should be $\Or(N_e^3)$. Numerically we observe that for the 1D examples up to $N_A=150$, the computational scaling is still in the pre-asymptotic regime. 

\subsection{2D lattice model}
In the previous section, we have validated the accuracy of ACP compared to both FD and DFPT. We also find that the efficiency of FD and DFPT can be comparable. Hence for the 2D model, we only compare the efficiency and accuracy of ACP with respect to DFPT.
Our first example is a periodic triangular lattice relaxed to the equilibrium position.
The distance between each atom and its nearest neighbor is set to be 1.2 a.u., and $\epsilon_0 = 0.05$. The nuclear charge $Z_I$ is set to 1 for all atoms, and $\sigma_I$ is set to be $0.24$.

\begin{figure}[!ht]
	\begin{minipage}[t]{0.5\linewidth}
		\centering
		\includegraphics[width=2.2in]{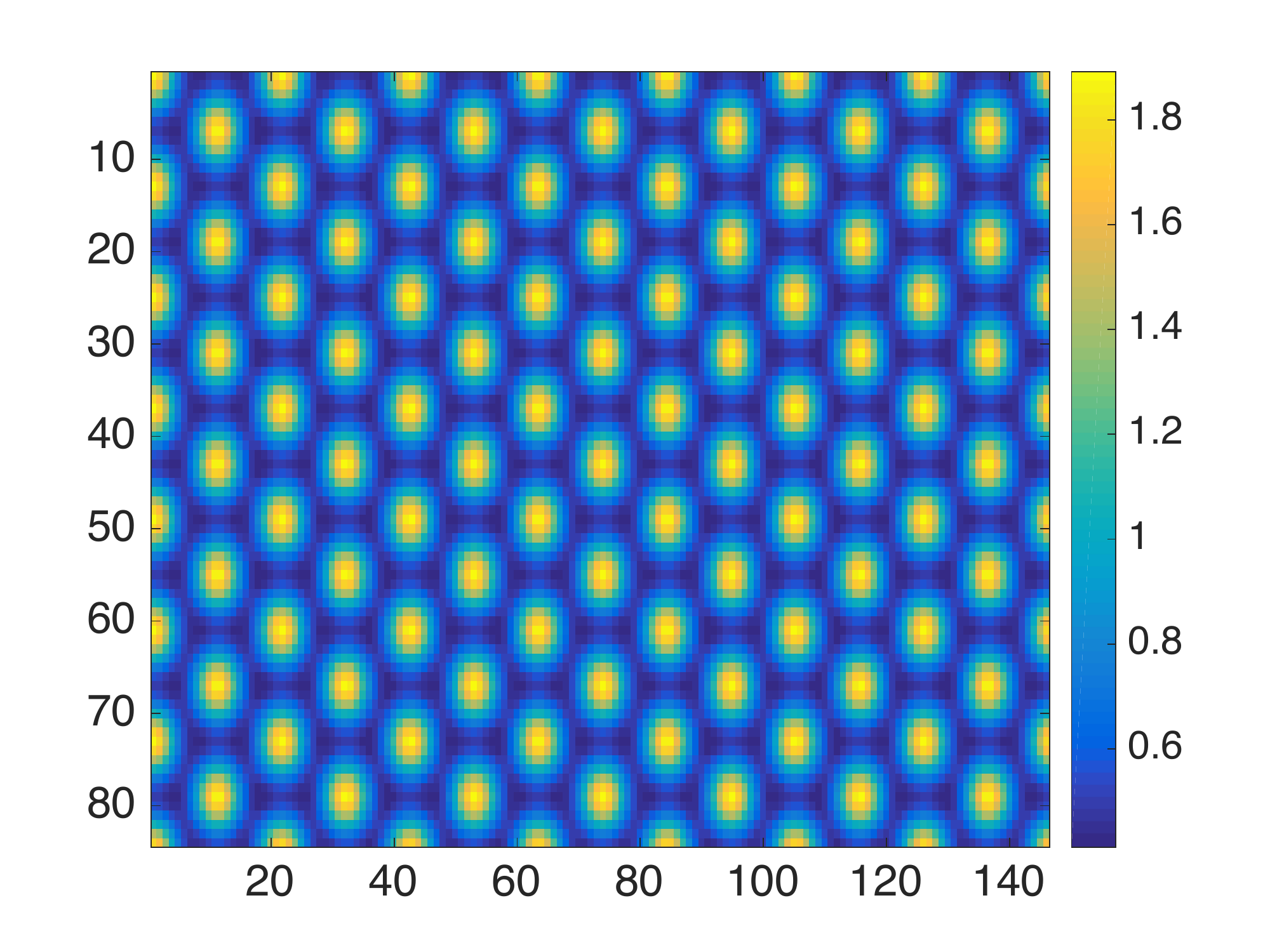}
		\subcaption{}
		\label{fig:rho77}
	\end{minipage}%
	\begin{minipage}[t]{0.5\linewidth}
		\centering
		\includegraphics[width=2.2in]{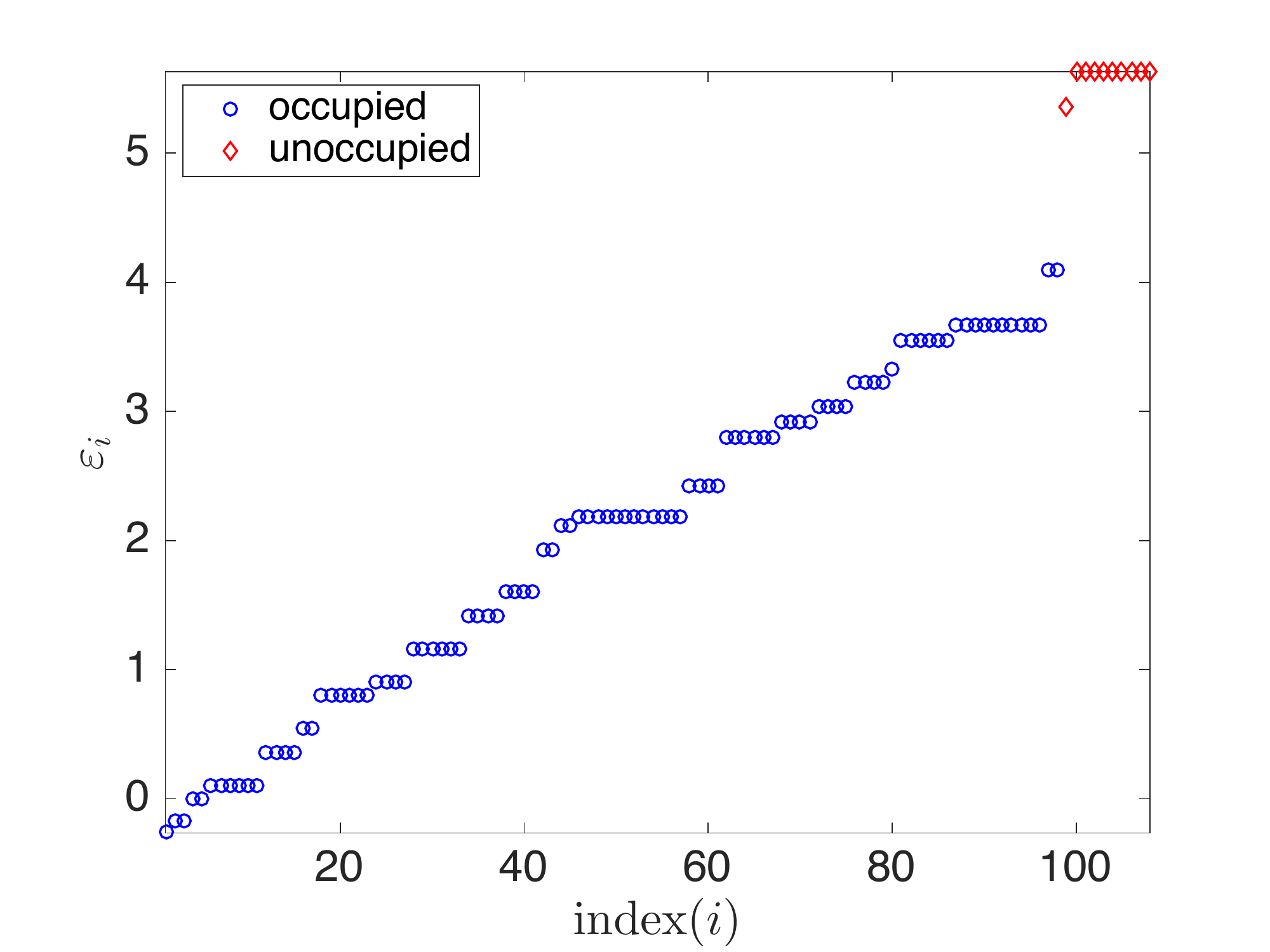}
		\subcaption{}
		\label{fig:gap77}
	\end{minipage}
	\caption{The electron density $\rho$ of the 98-atom insulating system (a), and the occupied and unoccupied eigenvalues (b).}
	\label{fig:2D}
\end{figure}

\begin{figure}[!ht]
	\begin{center}
		\includegraphics[width=0.50\textwidth]{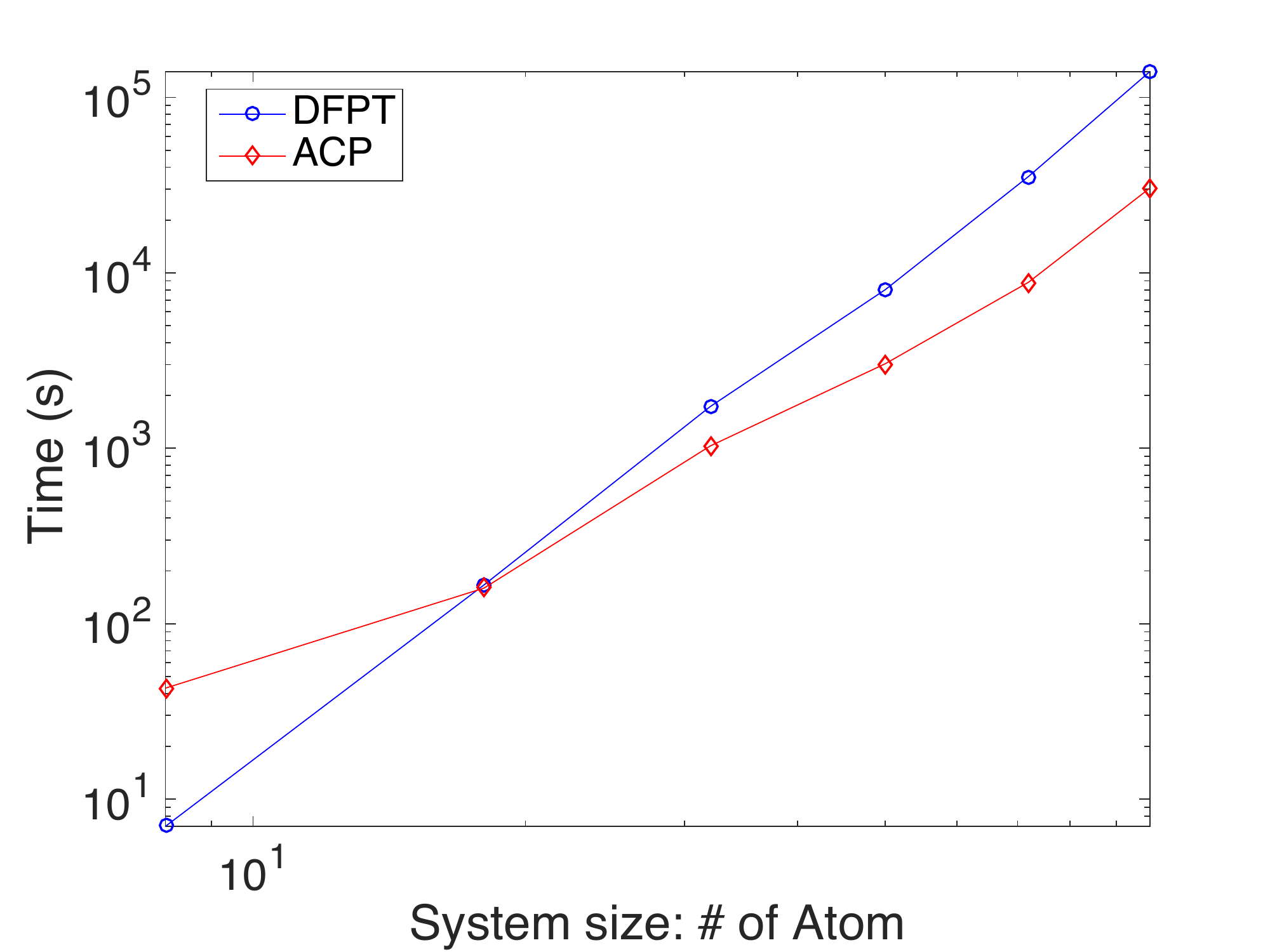}
	\end{center}
	\caption{Computational time. Comparison of ACP to ACP for the 2D periodic lattice model.}
	\label{fig:eff2D_1e-3}
\end{figure}

\begin{figure}[!ht]
	\begin{center}
		\includegraphics[width=0.50\textwidth]{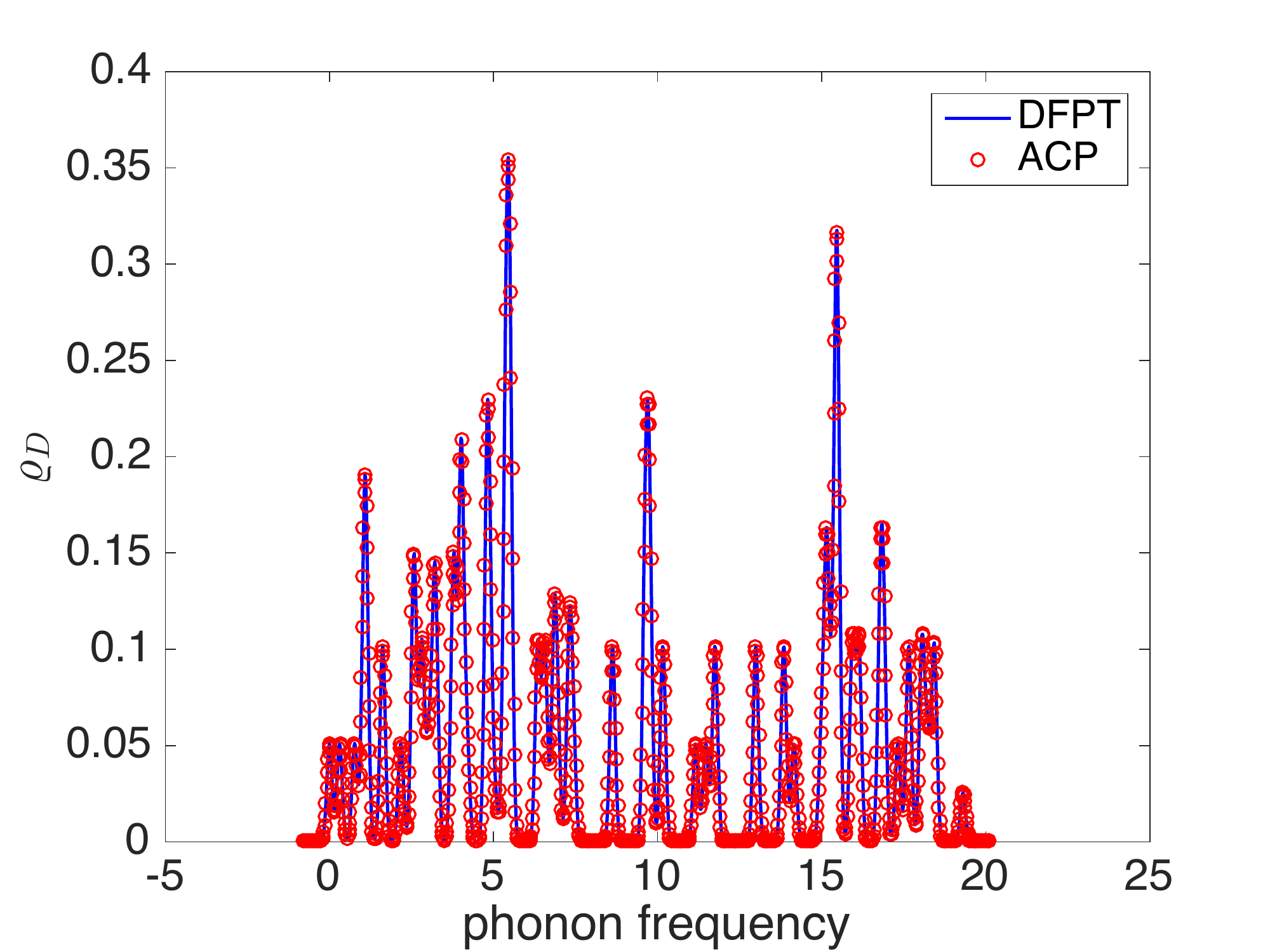}
	\end{center}
	\caption{Phonon spectrum for the 2D periodic lattice. System size $N_A = 98$. $\epsilon = 10^{-3}$, $N_\mu \approx 14 N_e$.}
	\label{fig:phononspec2D}
\end{figure}

For a system of size $N_A = 98$, the converged electron density $\rho$ as well as the 108 smallest eigenvalues associated with the Hamiltonian at the converged $\rho$ are shown in Fig.~\ref{fig:2D}. There is a finite HOMO-LUMO gap, $\varepsilon_g = \varepsilon_{99} - \varepsilon_{98} = 1.2637$, which suggests that the system is an insulator.
Fig.~\ref{fig:eff2D_1e-3} shows the computational time of ACP and DFPT as the system size grows from $2\times 2^2$ to $2\times 7^2$. We choose $\epsilon=10^{-3}$ in Alg.~\ref{alg:compressM}, and we find that this amounts to around $N_{\mu}= 14 N_e$ columns selected in the ACP procedure. We choose $N_c=30$, and hence the ACP formulation solves $420 N_e$ equations, compared to the $2N_A N_e = 2N_e^2$ equations needed for DFPT. We iterate Alg.~\ref{alg:overall} for two iterations.
We find that when the system size increases beyond $N_A=18$, the ACP formulation becomes more advantageous compared to DFPT. For  the largest system $N_A = 98$, ACP is 4.61 times faster than DFPT. Table \ref{tab:slope2D} measures the computational scaling from $N_A = 32$ to $N_A = 98$, which matches tightly with $\Or(N_e^3)$ and $\Or(N_e^4)$ theoretical scaling of ACP and DFPT, respectively.   Fig.~\ref{fig:phononspec2D} reports the phonon spectrum $\varrho_D$ for the system of size $N_A = 98$. Here the smear parameter $\sigma=0.08$. 



\begin{table}[!ht]
	\begin{center}
		\begin{tabular}{c|c} \hline
			Method & Slope   \\ \hline
			DFPT  & 3.9295  \\ \hline
			ACP  & 3.0249 \\ \hline
		\end{tabular}
	\end{center}
	\caption{Computational scaling measured from $N_A = 32$ to $N_A = 98$. }
	\label{tab:slope2D}
\end{table}


\subsection{2D model with random vacancies}

Our final example is a 2D triangular lattice with defects. We start from a periodic system with $N_A=72$ atoms, randomly remove three atoms, and then perform structural relaxation for $15$ steps. We terminate the structural relaxation before the system reaches its equilibrium position to obtain a disordered structure.

\begin{figure}[!ht]
	\begin{minipage}[t]{0.5\linewidth}
		\centering
		\includegraphics[width=2.2in]{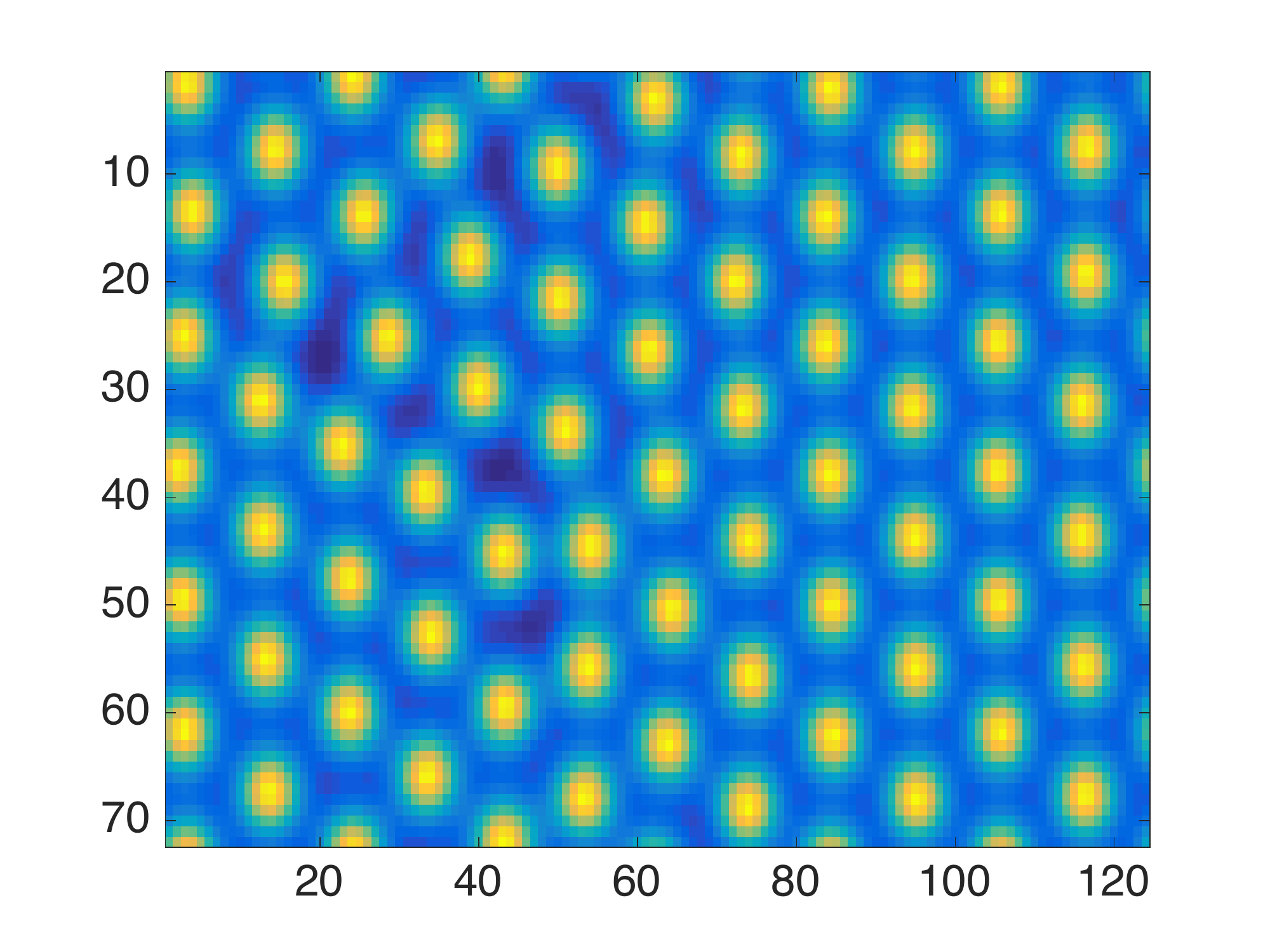}
		\subcaption{}
		\label{fig:rho66perturbed}
	\end{minipage}%
	\begin{minipage}[t]{0.5\linewidth}
		\centering
		\includegraphics[width=2.2in]{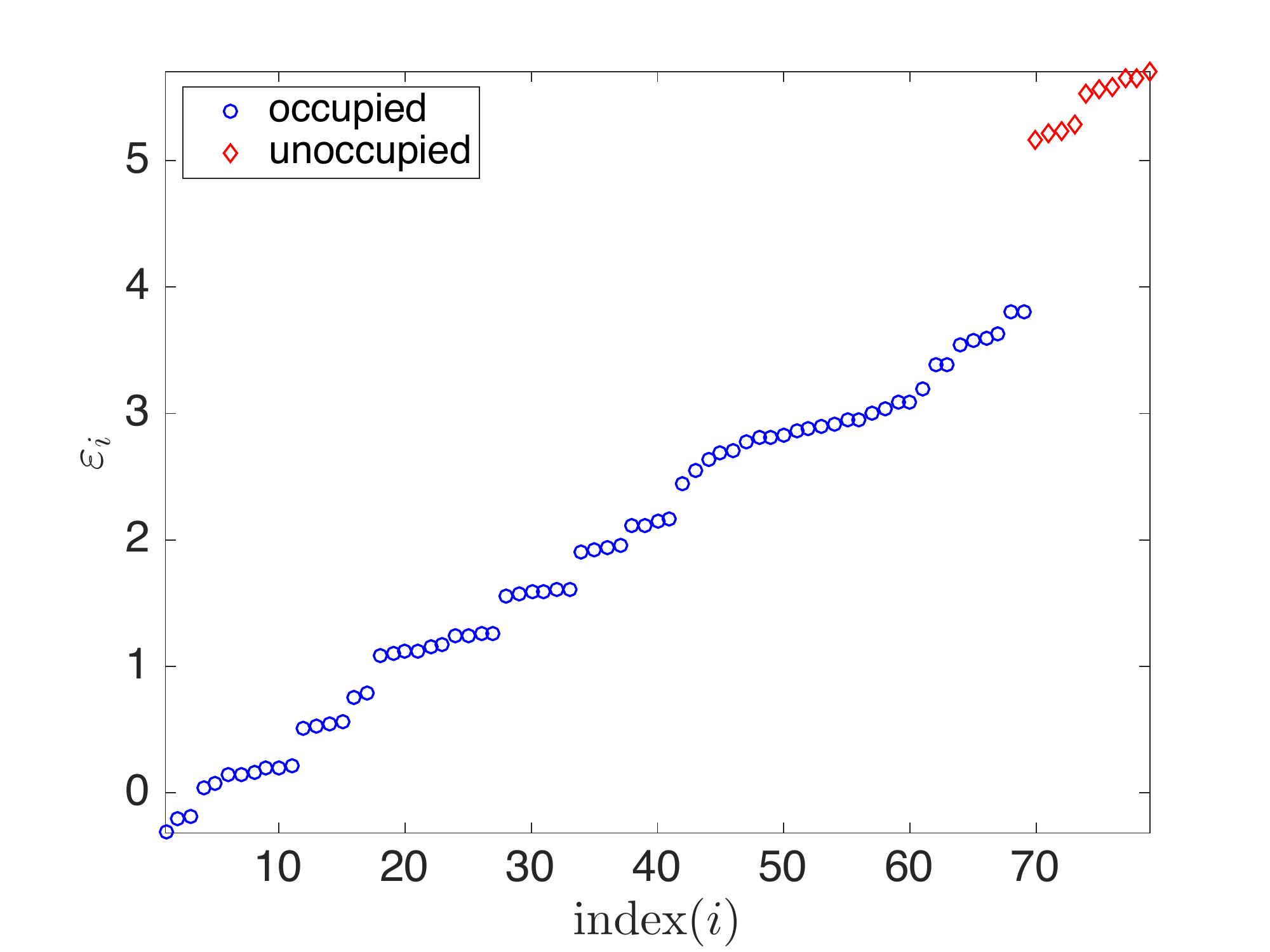}
		\subcaption{}
		\label{fig:gap66perturbed}
	\end{minipage}
	\caption{The electron density $\rho$ of the 2D system with defects (a), and the occupied and unoccupied eigenvalues (b).}
	\label{fig:2Dperturbed}
\end{figure}

\begin{figure}[!ht]
	\begin{center}
		\includegraphics[width=0.50\textwidth]{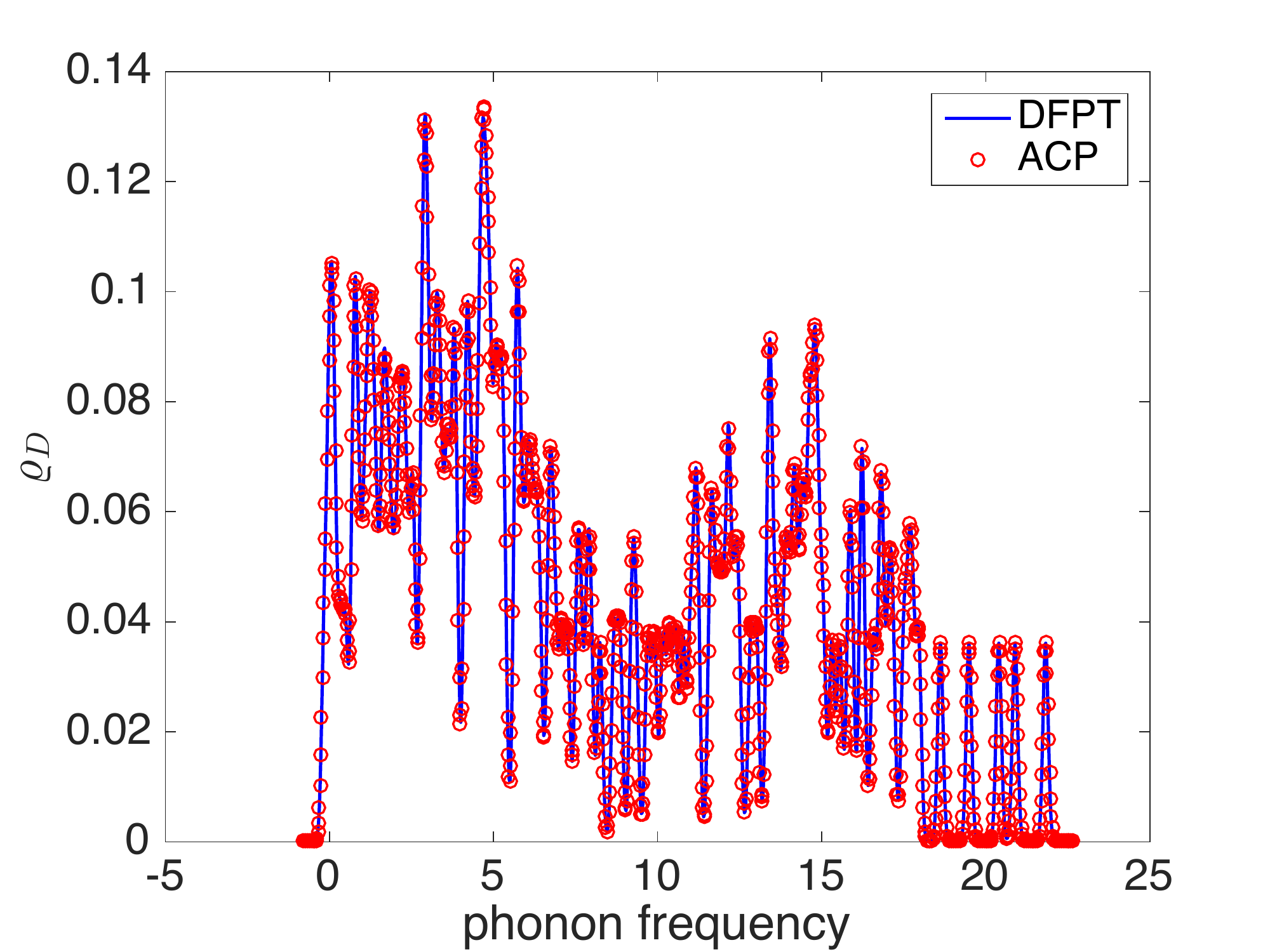}
	\end{center}
	\caption{Phonon spectrum for the 2D system with defects. System size $N_A = 69$, $\epsilon = 10^{-3}$, $N_\mu \approx 14 N_e$.}
	\label{fig:spec66perturbed}
\end{figure}

The converged electron density $\rho$ and the smallest 79 eigenvalues are shown in Fig.~\ref{fig:2Dperturbed}. For this system, there is a finite gap $\varepsilon_g = \varepsilon_{70} - \varepsilon_{69} = 1.3500$.
Fig.~\ref{fig:spec66perturbed} shows the  phonon spectrum computed from ACP and DFPT, plotted with the smear parameter $\sigma=0.08$. The computational time for DFPT is 53883 sec and that for ACP is 8741 sec, and the speedup factor for ACP is 6.16. 
We observe that for the disordered structure, DFPT requires more iterations to converge, while the number of iterations for ACP to converge can remain to be chosen to be 2. 
More specifically, compared to the periodic structure with $N_A=72$, the computational time for DFPT is 35501 sec, while that for ACP is 8854 sec. 

\section{Conclusion}\label{sec:conclusion}

We have introduced the adaptively compressed polarizability operator (ACP) formulation. To the extent of our knowledge, the ACP formulation reduces the computational complexity of phonon calculations from $\Or(N_e^4)$ to $\Or(N_e^3)$ for the first time.  
This is achieved by reducing the $\Or(N_e^2)$ equations in density functional perturbation theory (DFPT) to $\Or(N_e)$ equations with systematic control of  accuracy. Moreover, the accuracy of the ACP formulation depends weakly on the size of the gap, and hence can be applied to both insulator and semiconductor systems. Our numerical results for model problems indicate that the computational advantage of the ACP formulation can be clearly observed compared to DFPT and finite difference, even for systems of relatively small sizes. 

In the current work, for simplicity we have not included the nonlocal 
contribution in a pseudopotential framework. We will present the ACP formulation in the presence of the nonlocal pseudopotential, and its application for computing the phonon spectrum for real materials
in the near future. The availability of fast phonon calculations provides a
possible way to accelerate  structural relaxation optimization of large
scale molecules and solids. 
In this work we have restricted
ourselves to zero temperature calculations. We plan to extend the ACP
formulation to treat systems at finite temperature and hence metallic
systems.
We have used phonon calculation as an example to demonstrate the effectiveness of the compressed polarizability operator. The same strategy can be applied to applications of DFPT other than phonon calculations. We also plan to extend the ACP formulation to treat frequency-dependent polarizability operators that arise from many body perturbation theories in the future. 




\section*{Acknowledgments}
This work was partially supported by Laboratory Directed Research and
Development (LDRD) funding from Berkeley Lab, provided by the Director,
Office of Science, of the U.S. Department of Energy under Contract No.
DE-AC02-05CH11231 (L. L. and Z. Xu), by the Scientific Discovery through
Advanced Computing (SciDAC) program and the Center for Applied
Mathematics for Energy Research Applications (CAMERA) funded by U.S.
Department of Energy, Office of Science, Advanced Scientific Computing
Research and Basic Energy Sciences, and by the Alfred P. Sloan
fellowship (L. L.). L. Y. is partially supported by the National Science Foundation under grant DMS-0846501 and the DOE’s Advanced Scientific Computing Research program under grant DEFC02-13ER26134/DESC0009409. L. L. thanks useful discussion with Martin Head-Gordon.


\bibliographystyle{siam}
\bibliography{reference}

\end{document}